\documentclass[12pt]{amsart}
\usepackage{amssymb}
\usepackage[all]{xy}

\newcommand{\Q}{{\mathbb {Q}}}

\newcommand{\G}{{\bf{G}}}
\newcommand{\T}{{\bf{T}}}

\newcommand{\R}{{\mathbb{R}}}
\newcommand{\Z}{{\mathbb{Z}}}
\newcommand{\C}{{\mathbb{C}}}

\newcommand{\se}{{\bf{S}}}
\newcommand{\de}{{\bf{D}}}
\newcommand{\te}{{\bf{T}}}
\newcommand{\be}{{\bf{B}}}
\newcommand{\pe}{{\bf{P}}}

\newcommand{\N}{{\mathbb{N}}}

\newcommand{\Ad}{{\operatorname{Ad}}}

\newcommand{\Lie}{\operatorname{Lie}}

\newcommand{\OO}{{\mathcal O}}
\newcommand{\SSS}{{\mathcal{S}}}
\newcommand{\RR}{{\mathcal{R}}}
\newcommand{\LL}{{\mbox{\boldmath $\mathfrak{g}$}}}
\newcommand{\HH}{{\mbox{\boldmath $\mathfrak{h}$}}}
\newcommand{\bb}{{\mbox{\boldmath $\mathfrak{b}$}}}
\newcommand{\UU}{{\mbox{\boldmath $\mathfrak{u}$}}}
\newcommand{\LLL}{{\mbox{\boldmath $\mathfrak{g}$}(\mathcal{O})}}

\newcommand{\rank}{{\rm rank}}

\theoremstyle{plain}
\newtheorem{thm}{Theorem}[section]
\newtheorem{lem}[thm]{Lemma}
\newtheorem{prop}[thm]{Proposition}
\newtheorem{cor}[thm]{Corollary}

\theoremstyle{definition}
\newtheorem{remark}[thm]{Remark}

\newtheorem{defn}[thm]{Definition}
\def\bydefn{\stackrel{def}{=}}

\title[Decomposable forms and tori orbits]{Values of decomposable forms at $S$-integer points and 
tori orbits on homogeneous spaces}

\author{George Tomanov}
\address{Institut Girard Desargues, Universit\'e Claude Bernard - Lyon
I, B\^atiment de Math\'ematiques, 43, Bld.
 du 11 Novembre 1918,
69622 Villeurbanne Cedex, France {\tt tomanov@igd.univ-lyon1.fr}}

\begin{document}

\maketitle

\begin{abstract}
Let $\G$ be a reductive algebraic group defined over a number field
$K$ and let $\SSS$ be a finite set of non-equivalent valuations of
$K$ containing the archimedean ones. Let $G = \prod_{v \in
\SSS}\G(K_v)$ and $\Gamma$ be an $\SSS$-arithmetic subgroup of $G$.
Let $\RR \subset \SSS$ and $T_\RR = \prod_{v \in \RR}T_v$ where
$T_v$ is a sub-torus of $\G(K_v)$ containing a maximal $K_v$-split
torus. We prove that if $G/\Gamma$ admits a closed $T_\RR$-orbit
then $\RR = \SSS$ or $\RR$ is a singleton. In addition, the closed
$T_\RR$-orbits are always "standard"; this generalizes the result of
\cite{Tomanov-Weiss}. When $\# \SSS
> 1$ it turns out that for $\RR = \SSS$ there are no divergent orbits and
for $\# \RR = 1$ all closed orbits are divergent. %Our study of the
%orbits of maximal tori is mainly motivated by problems from number
%theory.
As an application, we prove that if a collection of decomposable
homogeneous forms $f_v \in K_v[x_1, \ldots, x_n], v \in \SSS,$ takes
discrete values at $\OO^n$, where $\OO$ is the ring of
$\SSS$-integers of $K$, then there exists an homogeneous form $g \in
\OO[x_1, \ldots, x_n]$ such that $f_v = \alpha_v g$, $\alpha_v \in
K_v^*$, for all $v \in \SSS$.
\end{abstract}

\section{Introduction} \label{Introduction}

Let $\G$ be a reductive algebraic group defined over a number field
$K$ and let $\SSS$ be a finite set of (normalized) valuations of $K$
containing all archimedean ones. If $v \in \SSS$ we set $G_v =
\G(K_v)$, where $K_v$ is the completion of $K$ with respect to $v$.
Every $G_v$ is a locally compact group with a topology induced by
the topology of $K_v$. Let $G = \prod_{v \in \SSS}G_v$. The group of
$K$-rational points $\G(K)$ is identified with its diagonal
imbedding in $G$. We denote by $\Gamma$ an $\SSS$-arithmetic
subgroup of $G$, that is, $\Gamma$ is a subgroup of $G$ such that
$\Gamma \cap \G(\OO)$ has finite index in both $\Gamma$ and
$\G(\OO)$, where $\OO$ is the ring of $\SSS$-integers of $K$. We fix
a maximal $K$-split torus $\de$ of $\G$ and, for every $v \in \SSS$,
we fix a $K_v$-torus $\T_v$ of $\G$ such that $\T_v$ contains both
$\de$ and a maximal $K_v$-split torus of $\mathbf{G}$. Let
$\mathcal{R}$ be a  non-empty  subset of $\SSS$. Recall that the
$\RR$-{\it rank} of $\G$ (or $G$) is $\rank_\RR \G \bydefn \sum_{v
\in \RR} \rank_{K_v} \G$. (If $F$ is a field containing $K$ then
$\rank_{F} \G$ is by definition the dimension of any maximal
$F$-split torus of $\G$.) We set $T_\RR = \prod_{v \in \RR}T_v$ and
$D_\RR = \prod_{v \in \RR}D_v$, where $T_v = \te_v(K_v)$ and $D_v =
\de(K_v)$. Then $T_\RR$ is a {\it torus of maximal
$\mathcal{R}$-rank} and it acts on $G/\Gamma$ by left translations
$$
t\pi(g) = \pi(tg),
$$
where $\pi: G \rightarrow G/\Gamma$ is the quotient map. An orbit
$T_\RR\pi(g)$ is called {\it divergent} if the orbit map $t
\rightarrow t\pi(g)$ is proper, i.e. if $\{t_i\pi(g)\}$ leaves
compacts of $G/\Gamma$ whenever $\{t_i\}$ leaves compacts of
$T_\RR$. In particular, the divergent orbits are closed.

We prove the following:

\begin{thm}
\label{thm1} Let $\rank_\RR \G > 0$ and $g \in G$.
\begin{enumerate}
\item[(a)] The orbit $T_\RR \pi(g)$ is closed if and only if $\RR$ is a singleton or $\RR =
\SSS$, and there exists a $K$-torus $\mathbf{L}$ of $\G$ such that
$$g^{-1}T_\RR g = C L_\RR,$$
where $C$ is a compact group and $L_\RR = \prod_{v \in \RR}
\mathbf{L}(K_v)$;
%and $T_\RR \pi(g)$ is divergent, or $\RR = \SSS$ and there exists a
%$K$-torus $\mathbf{L}$ of $\G$ such that
%$$g^{-1}T_\RR g = C L,$$
%where $C$ is a compact subgroup and $L = \mathbf{L}(K_\SSS)$; %with
\item[(b)] The orbit $T_\RR \pi(g)$ is divergent if and only if the following conditions are satisfied:
$\RR$ is a singleton equal to $v$, $\rank_{K_v} \G = \rank_K \G$ and
$$g \in \mathcal{Z}_{G}(D_v)\G(K),$$
where $D_v$ is identified with its natural projection in $G$ and
$\mathcal{Z}_{G}(D_v)$ is the centralizer of $D_v$ in $G$.
\end{enumerate}
\end{thm}

Theorem \ref{thm1} generalizes the following result by B.Weiss and
the author, the second part of which has been earlier proved (though
unpublished) by G.Margulis for $\G = \mathbf{{SL}_n}$
endowed with the standard $\Q$-structure %$\Gamma = \mathrm{SL}_n(\Z)$
(cf. \cite[Appendix]{Tomanov-Weiss}).

\begin{thm}
\label{ToW} \mbox{\rm{(}}\cite[Theorem
1.1]{Tomanov-Weiss}\mbox{\rm{)}} Let $\G$ be a reductive
$\Q$-algebraic group, $\mathbf{T}$ an
 $\R$-torus containing   a maximal $\R$-split torus,
$T=\mathbf{T}(\R)$ and let $x \in G$. Then:
\begin{itemize}
\item
$T \pi(x)$ is a closed orbit if and only if  $x^{-1}\T{x}$ is a
product of a $\Q$-subtorus and an  $\R$-anisotropic $\R$-subtorus;
\item
$T\pi(x)$ is a divergent orbit if and only if the maximal $\R$-split
subtorus of $x^{-1}\T{x}$ is defined over $\Q$ and $\Q$-split.
\end{itemize}
\end{thm}

When $\#\RR > 1$, Theorem \ref{thm1} implies a specific phenomenon:

\begin{cor}
\label{S-div} If $\#\SSS > 1$ and $T_\RR \pi(g)$ is a closed orbit
then either $\RR = \SSS$ and $T_\RR \pi(g)$ is never divergent, or
$\RR$ is a singleton and $T_\RR \pi(g)$ is always divergent.
\end{cor}

An orbit $T_\RR \pi(g)$ is called {\it locally divergent} if
$T_{v}\pi(g)$ is divergent for every $v \in \RR$. Theorem \ref{thm1}
will be deduced from the next theorem about the locally divergent
orbits.

\begin{thm}
\label{thm2} Let $\mathrm{rank}_{\RR}(\G) > 0$. %With the above notation, %from Theorem \ref{thm1},
Then the orbit $T_\RR \pi(g)$ is closed and locally divergent if and
only if the following conditions are fulfilled:
\begin{enumerate}
\item[(i)] $\RR = \SSS$ or $\RR$ is a singleton;
\item[(ii)]  %$\mathrm{rank}_{\RR}(\G) = |\RR| \mathrm{rank}_{K}(\G)$
$\mathrm{rank}_{\RR}(\G) = \#\RR \ \mathrm{rank}_{K}(\G)$;
%where $|\RR|$ is the cardinality of $\RR$;
%$\mathrm{rank}_{K_v}(\G) = \mathrm{rank}_{K}(\G)$ for every $v \in
%\RR$;
%, i.e. $T_v/D_v$ is a compact group;
%There exists a $v_0 \in \SSS$ such that $T_{v_{0}}x$ is
%divergent;
\item[(iii)]
$g \in \mathcal{N}_{G}(D_\RR)\G(K)$, where $\mathcal{N}_{G}(D_\RR)$
is the normalizer of $D_\RR$ in $G$.
\end{enumerate}
\end{thm}

When $\#\RR = 1$ we can replace the normalizer
$\mathcal{N}_{G}(D_\RR)$ in the formulation of Theorem \ref{thm2}
(iii) by the centralizer $\mathcal{Z}_{G}(D_\RR)$.  This is not
possible when $\RR = \SSS$ (see 6.2 (b)).

As a consequence of Theorem \ref{thm2}, one can easily see that the
locally divergent $T_\RR$-orbits are also all "standard":

\begin{cor}
\label{thm3} Let $g \in G$. The orbit $T_\RR \pi(g)$ is locally
divergent if and only if $$\mathrm{rank}_{\RR}(\G) = \#\RR \
\mathrm{rank}_{K}(\G)$$ and $$g \in \bigcap_{v \in
\RR}\mathcal{Z}_{G}(D_{v})\G(K).$$
\end{cor}

We also get the following result:
\begin{cor}
\label{thm3'}% Let $g = (g_v)_{v \in \SSS} \in G$. The following
%assertions hold:
\begin{enumerate}
\item[(a)] If $\mathrm{rank}_{\RR}(\G) > \#\RR \ \mathrm{rank}_{K}(\G)$ then there are no
locally divergent orbits for $T_\RR$;
\item[(b)] Let $\G$ be semisimple, $\#\RR > 1$ and \ $\mathrm{rank}_{\RR}(\G) = \#\RR \ \mathrm{rank}_{K}(\G)
> 0$. %$\rank_{K_{v}} \G  = \rank_{K} \G > 0$ for all $v \in \RR$.
Then there exist locally divergent but non-closed orbits for
$T_\RR$.
\end{enumerate}
\end{cor}

We apply Theorem \ref{thm1} to obtain a characterization of the
rational decomposable homogeneous forms in terms of their values at
the integer points. Such forms appear in a very natural way in both
the algebraic number theory and the Diophantine approximation of
numbers in connection with the notable Littlewood conjecture. (See,
\cite[ch.2]{BoSha} and \cite[\S 2]{survey}, respectively.)

We will first formulate our result in technically simpler particular
cases. Given a commutative ring $R$, we denote by $R[\ \vec{x}\ ]$
the ring of polynomials with coefficients from $R$ in $n$ variables
$\vec{x} = (x_1,\ldots,x_n)$.

\begin{thm}
\label{inverse} Let $f(\vec{x}) = l_{1}(\vec{x})\ldots
l_{m}(\vec{x})$, where $l_{1}(\vec{x}),\ldots, l_{m}(\vec{x}) \in
\R[\ \vec{x}\ ]$ are real linear forms. Suppose that
$l_{1}(\vec{x}),\ldots, l_{m}(\vec{x})$ are linearly independent
over $\R$ and that the set $f(\Z^n)$ is discrete in $\R$. Then
$f(\vec{x}) = \alpha g(\vec{x})$, where $g(\vec{x}) \in
\Z[\,\vec{x}\,]$ and $\alpha \in \R^*$.
\end{thm}

The hypotheses that the form $f(\vec{x})$ is decomposable and
$l_{1}(\vec{x}),\ldots,$ $ l_{m}(\vec{x})$ are linearly independent
over $\R$ are essential. (See \S 7 for simple examples.) It is easy
to prove (see \cite[ch.2, Theorem 2]{BoSha}) that the form
$g(\vec{x})$ in the formulation of the theorem is a constant
multiple of a product of forms of the type $\mathrm{N}_{K/\Q}(x_1 +
x_2\mu_2 + \ldots + x_n\mu_n)$, where $\mu_2, \ldots, \mu_n$ are
algebraic numbers linearly generating a totally real number field
$K$ of degree $n$ and $\mathrm{N}_{K/\Q}$ is the algebraic norm of
$K$.

If $f$ is a decomposable homogeneous form with complex coefficients
and we are considering the values of $f$ at the Gaussian integer
vectors, % $\Z[\,i\,]^n$
we get:

\begin{thm}
\label{C-inverse} Let $f(\vec{x}) = l_{1}(\vec{x})\ldots
l_{m}(\vec{x})$, where $l_{1}(\vec{x}),\ldots, l_{m}(\vec{x}) \in
\C[\ \vec{x}\ ]$ are complex linear forms. Suppose that
$l_{1}(\vec{x}),\ldots, l_{m}(\vec{x})$ are linearly independent
over $\C$ and that the set $f(\Z[\,i\,]^n)$ is discrete in $\C$.
Then $f(\vec{x}) = \alpha g(\vec{x})$, where $g(\vec{x}) \in
\Z[\,i\,][\,\vec{x}\,]$ and $\alpha \in \C^*$.
\end{thm}

Let $K$, $\SSS$ and $\OO$ be as in the formulation of Theorem
\ref{thm1}. For every $v \in \SSS$, let $f_v = l_{1}^{(v)} \ldots
l_{m}^{(v)} \in K_v[\,\vec{x}\,]$, where $l_{1}^{(v)}, \ldots,
l_{m}^{(v)}$ are linearly independent over $K_v$ linear forms in
$K_v[\,\vec{x}\,]$. Denote by $K_\SSS$ the direct product of the
topological fields $K_v$, $v \in \SSS$.
%Then $K_\SSS$ is a topological ring.
Both Theorems \ref{inverse} and
\ref{C-inverse} are particular cases for $K = \Q$ and $K =
\Q(\,i\,)$, respectively, of the next general theorem:

\begin{thm}
\label{S-inverse} With the above notation, assume that
$\{(f_v(\vec{z}))_{v \in \SSS} \in K_\SSS | \vec{z} \in \OO^n \}$ is
a discrete subset of $K_{\SSS}$. Then there exist an homogeneous
form $g$ with coefficients from $\OO$ and an element $(\alpha_v)_{v
\in \SSS} \in K_\SSS^*$ such that $f_v = \alpha_v g$ for all $v \in
\SSS$.
\end{thm}

In connection with Theorem \ref{S-inverse} it seems natural to
formulate the following conjecture which generalizes a well known
conjecture for the real forms $f$:

\vspace{.3cm}

$\mathbf{Conjecture.}$ Let $f_v, v \in \SSS$, be as in the
formulation of Theorem \ref{S-inverse} with $n = m$ and $\#
\SSS.n > 2$. Additionally, assume that %$m \geq 3$ and that
there exists a neighborhood $W$ of $0$ in $K_\SSS$ such that
$(f_v(\vec{z}))_{v \in \SSS} \notin W$ for every $\vec{z} \in \OO^n,
\vec{z} \neq 0$. Then there exist an homogeneous form $g$ with
coefficients from $\OO$ and an element $(\alpha_v)_{v \in \SSS} \in
K_\SSS^*$ such that $f_v = \alpha_v g$ for all $v \in \SSS$.

\vspace{.3cm}

Using the $\SSS$-adic version of Malher's criterion (see Theorem
\ref{mahler} below), it is easy to see that the above conjecture can
be reformulated in terms of Theorem \ref{thm1} as follows: If $\G =
\mathbf{SL}_n$ and $\rank_\SSS \G > 1$ then $T_\SSS \pi(g)$ is
compact whenever $T_\SSS \pi(g)$ is relatively compact. In the case
$n = 3$ and $K = \Q$ the conjecture implies (cf.\cite[\S 2]{survey})
the Littlewood conjecture which states that
$$
\liminf_{n \rightarrow \infty} n\langle n \alpha \rangle \langle
n\beta \rangle = 0
$$
for all $\alpha, \beta \in \R$, where $\langle x \rangle$ denotes
the distance from $x$ to $\Z$. In \cite{EiKaLi}, using the dynamical
approach, M.Einsiedler, A.Katok and E.Lindenstrauss proved that the
Littlewood conjecture fails at most on a set of Hausdorff dimension
zero. Similar results in the $p$-adic setting have recently appeared
in the M.Einsiedler and D.Kleinbock paper \cite{EiKl}.

 \vspace{.6cm}

The paper is organized as follows. The notation and the terminology
are introduced in \S 2. Our starting point is the paper
\cite{Tomanov-Weiss}. In \S 3, using \cite{Tomanov-Weiss}, we prove
an $\SSS$-adic compactness criterium in terms of intersections of
so-called quasiballs with horospherical subsets. In \S 4 we prove
Proposition \ref{K-torus} which plays a crucial role in revealing
the dichotomy in Corollary \ref{S-div}. In \S 5 we describe the
locally divergent orbits in terms of minimal parabolic $K$-algebras.
In order to do this, we have to apply more intrinsic arguments than
in \cite[\S 5]{Tomanov-Weiss} for the proof of a similar result. For
instance, the Galois type arguments are replaced by Proposition
\ref{prop7'} from the algebraic group theory. Theorems \ref{thm1},
\ref{thm2} and their corollaries are proved in \S 6. The proof of
Theorem \ref{S-inverse} is given in \S 7.

The author is grateful to Manfred Einsiedler, Dima Kleinbock,
Gregory Margulis and Barak Weiss for the useful discussions and to
the Max Planck Institut f\"ur Mathematik, where the main part of
this work was accomplished, for its hospitality.

\section{Preliminaries: notation and basic concepts} \label{section:
prelims}

\subsection{Numbers}

As usual $\C$, $\R$, $\Q$ and $\Z$ denote the complex, real,
rational and integer numbers, respectively.

In this paper ${K}$ denotes a number field, that is, a finite
extension of $\Q$. All valuations of $K$ which we consider are
supposed to be \textit{normalized} (see \cite[ch.2, \S7]{CF}) and,
therefore, pairwise non-equivalent. If $v$ is a valuation of $K$
then ${K}_v$ is the completion of $K$ with respect to $v$ and $|\ .\
|_v$ is the
corresponding norm on $K_v$. %If $v$ is
%archimedean (respectively, non-archimedean) then ${K}_v$ is
%isomorphic to  $\R$ or  $\C$ (respectively, to a finite extension of
%the field of $p$-adic numbers $\Q_p$).
If $v$ is non-archimedean then ${\OO}_v = \{x \in K_v : |\
 x\ |_v \leq 1 \}$ is the ring of integers of ${K}_v$.

We fix a finite set $\mathcal{S}$ of valuations of ${K}$ containing
all archimedean valuations of $K$. The latter set is denoted by
$\mathcal{S}_\infty$ or, simply, $\infty$, if this does not lead to
confusion. We also put $\mathcal{S}_f = \mathcal{S} \setminus
\mathcal{S}_\infty$.

We denote by $\mathcal{O}$ the ring of $\mathcal{S}$-integers of
$K$, i.e., $\mathcal{O} = K \bigcap (\bigcap_{v \notin
\SSS}\mathcal{O}_v)$.

For any non-empty subset $\mathcal{R}$ of $\mathcal{S}$,
${K}_\mathcal{R} \bydefn \prod_{v \in \mathcal{R}}{K}_v$ is a direct
product of locally compact fields. Note that ${K}_\mathcal{R}$ is a
topological ring and that the diagonal imbedding of $K$ in
$K_\mathcal{R}$ is dense. As usual, we denote by $K_\mathcal{R}^*$
the multiplicative group of all invertible elements in the ring
$K_\mathcal{R}$.

\subsection{Norms}
Let $\mathbf{V}$ be a finite dimensional vector space defined over
$K$. For every $\RR \subset \SSS$  (respectively $v \in \SSS$) we
write $V_\RR$ for $\mathbf{V}(K_\RR)$ (respectively, $V_v$ for
$\mathbf{V}(K_v)$). Fixing a basis of $K$-rational vectors $e_1,
\ldots, e_n$, for every $K$-algebra $A$, we identify $\mathbf{V}(A)$
with $A^n$. For every $v \in \SSS$ we define a {\it normalized} norm
$\|\cdot\|_v$ on $V_v$ as follows. If $v$ is real (respectively,
complex) then $\|\cdot\|_v$ is the standard norm on $\R^n$
(respectively, the square of the standard norm on $\C^n$). If $v$ is
non-archimedean, then $\|\cdot\|_v$ is defined by $\|\mathbf{x}\|_v
= \max_i|x_i|_v$, where $(x_1,\dots,x_n)$ are the coordinates of the
vector $\mathbf{x} \in V_v$ with respect to the bases $e_1, \ldots ,
e_n$.

For $\textbf{x} = (\textbf{x}^{(v)})_{v \in \SSS}$ in $V_\RR$ we
define the norm of $\mathbf{x}$ as
$$\|\mathbf{x}\|_\RR = \max_{v \in \RR} \|\mathbf{x}^{(v)}\|_v.$$

Also, if $\RR = \SSS$ we define the {\it content} of $\textbf{x}$ as
$$\mathbf{c}_\SSS(\textbf{x}) = \prod_{v \in \SSS} \|\mathbf{x}^{(v)}\|_v.$$
Since all our norms are normalized and $\prod_{v \in \SSS}|\xi|_v =
1$ for every $\xi \in \OO^*$ \cite[ch.2, Theorem 12.1]{CF}, we have
that
\begin{equation}
\label{Artin} \mathbf{c}_\SSS(\textbf{x}) =
\mathbf{c}_\SSS(\xi\textbf{x}), \forall \xi \in \OO^*.
\end{equation}

By a {\it pseudoball} in $V_\SSS$ of radius $r > 0$ centered at $0$
we mean the set $\mathcal{B}_\SSS(r) =
 \{\mathbf{x} \in V_\SSS |\mathbf{c}_\SSS(\textbf{x}) < r \}$. %We will write $\mathcal{B}_{\infty}(r)$ to
%denote the pseudoball of radius $r$ in $V_{\infty}$.
We preserve the notation $B_\SSS (r)$ to denote the usual ball in
$V_\SSS$ of radius $r$  centered at $0$ with respect to the norm $\|
. \|_\SSS$.

\subsection{$K$-algebraic groups and their Lie algebras} %Concerning the standard facts about the theory of linear
%algebraic groups we refer the reader to \cite{Borel}.
We use boldface upper case letters to denote the algebraic groups
and boldface lower case Gothic letters to denote their Lie algebras.

In this paper $\G$ is a reductive algebraic group defined over $K$.
Recall that the Lie algebra $\LL$ of $\G$ is equipped with a
$K$-structure compatible with the $K$-structure of $\G$
\cite[Theorem 3.4]{Borel}. An algebraic subgroup of $\G$ defined
over $K$ is called shortly $K$-{\it subgroup}.

Given $\RR \subset \SSS$ and  a $K$-subgroup $\mathbf{H}$ of $\G$,
we usually denote $H_\RR \bydefn \mathbf{H}(K_\RR)$ and
$\mathfrak{h}_\RR \bydefn \mbox{\boldmath $\mathfrak{h}$}(K_\RR)$.
The group $H_\RR$ (respectively, its Lie algebra $\mathfrak{h}_\RR$)
is identified with the direct product $\prod_{v \in \RR} H_v$
(respectively, $\prod_{v \in \RR}\mathfrak{h}_v$), where $H_v
\bydefn \mathbf{H}(K_{v})$ (respectively, $\mathfrak{h}_v \bydefn
\mbox{\boldmath $\mathfrak{h}$}(K_v)$). But if $\RR = \SSS$ and this
does not lead to confusion we prefer the simpler notation $H$
(respectively, $\mathfrak{h}$) for $H_\SSS$ (respectively,
$\mathfrak{h}_\SSS$).

We will use the notation $\mathrm{pr}_\RR$ to denote both the
natural projections $ G \to G_\RR$ and $\mathfrak{g} \to
\mathfrak{g}_\RR$. (The exact use of $\mathrm{pr}_\infty$ will
follow from the context.)

On every $G_v$ we have a {\it Zariski topology} induced by the
Zariski topology on $\bf G$ and a \textit{Hausdorff topology}
induced by the locally compact topology on $K_v$. The formal product
of the Zariski (respectively, Hausdorff) topologies on $G_v$, $v \in
\RR$, is the Zariski (respectively, Hausdorff) topology on $G_\RR$.
In order to distinguish the two topologies, all topological notions
connected with the first one will be used with the prefix "Zariski".

An element $g = (g_v)_{v \in \RR} \in G_\RR$ is called {\it
unipotent} (respectively, {\it semisimple}) if each $v$-component
$g_v$ of $g$ is unipotent (respectively, semisimple). A subgroup $U$
of $ G_\RR$ is called unipotent if it consists of unipotent
elements. A subalgebra  $\mathfrak{u}$ of $\mathfrak{g}_\RR$ is {\it
unipotent} if it corresponds to a Zariski closed unipotent subgroup
$U$ of $ G_\RR$, i.e. if there exists a subgroup $U \subset G_\RR$
such that $U = \prod_{v \in \RR}U_v$, each $U_v$ is Zariski closed
in $G_v$, and $\mathfrak{u} = \prod_{v \in \RR}\mathfrak{u}_v$ where
$\mathfrak{u}_v$ is the Lie algebra of $U_v$.

If $\mathbf{P}$ is a parabolic $K$-subgroup of $\G$ then
${R}_u(\mathbf{P})$ denotes the unipotent radical of $\mathbf{P}$.
The {\it unipotent radical of the Lie algebra of $\pe$} is by
definition the Lie algebra of ${R}_u(\mathbf{P})$.

If $H$ is a subgroup of  $G$ then $\mathcal{N}_{G}{(H)}$
(respectively, $\mathcal{Z}_G(H)$) denotes the normalizer
(respectively, the centralizer) of $H$ in $G$.

For any non-empty $\RR \subset \SSS$ the adjoint representation
$\mathrm{Ad}_\RR: G_\RR \rightarrow \mathrm{GL}(\mathfrak{g}_\RR)$,
where $\mathrm{GL}(\mathfrak{g}_\RR) = \prod_{v \in
\RR}\mathrm{GL}(\mathfrak{g}_v)$, is the direct product of the
adjoint representations $\mathrm{Ad}_v: G_v \rightarrow
\mathrm{GL}(\mathfrak{g}_v)$, $v \in \RR$. We will use the notation
$\Ad$ (respectively, $\Ad_\infty$) when $\RR = \SSS$ (respectively,
$\RR = \SSS_\infty$).

\subsection{$\SSS$-arithmetic subgroups}
Recall that  $\Gamma$ is  an $\SSS$-{\it arithmetic} subgroup of
$G$, i.e., $\Gamma \cap \G(\OO)$ has finite index in both $\Gamma$
and $\G(\OO)$. We assume that $\G$ is imbedded in $\mathbf{SL}_n$ in
such a way that $\G(\OO) = \mathbf{SL}_n(\OO) \cap \G$ and $\LL(\OO)
=$ ${\mbox{\boldmath $\mathfrak{sl}$}}$$_n(\OO) \cap \LL$. In
particular, $\LL(\OO)$ {\it is invariant under the adjoint action
of} $\G(\OO)$. Let $\Gamma'$ be a subgroup of finite index in
$\Gamma$ and let $\phi: G/\Gamma' \rightarrow G/\Gamma$ be the
natural map. Since $\phi$ is a proper map it is easy to see that
Theorems \ref{thm1}, \ref{thm2} and their corollaries are valid for
$\Gamma$ if and only if they are valid for $\Gamma'$. Therefore, we
may suppose {\it without} loss of generality that $\Gamma =
\G(\OO)$.

Let $\pi: G \rightarrow G/\Gamma$ be the natural projection. For
every $x \in G/\Gamma$ we introduce the following notation. If $x =
\pi(g)$, $g \in G$, we denote
$$\mathfrak{g}_x=\Ad(g)\LLL.$$
Since $\LLL$ is $\Ad (\Gamma)$-invariant, $\mathfrak{g}_x$ does not
depend on the choice of the element $g$.

\section{Compactness criteria in $\SSS$-adic setting}

%such that
%$\rank_{\infty} \mathbf{G}_i > 0$ for every non-trivial $K$-simple
%component of $\G$.

\subsection{$\SSS$-adic  Mahler's criterion}

Let $G = \mathrm{SL}_n(K_\SSS)$, $\Gamma = \mathrm{SL}_n(\OO)$ and
$\pi: G \rightarrow G/\Gamma$ be the natural projection. The group
$G$ is acting naturally on $K_\SSS^n$ and $\Gamma$ is the stabilizer
of $\OO^n$ in $G$. If $r > 0$ then $B_\SSS(r)$ (resp.,
$\mathcal{B}_\SSS(r)$) is the ball (resp. pseudoball) in $K_\SSS^n$
centered in $0$ and with radius $r$ (see \S 2.3).

We have

\begin{thm}
\label{mahler} $\mathrm{(Mahler's \  criterion)}$ With the above
notation, given a subset $M \subset G$ the following conditions are
equivalent:
\begin{enumerate}
\item[(i)] $\pi(M)$ is relatively compact in
$G/\Gamma$;
\item[(ii)] There exists $r > 0$ such that $g\OO^n \cap
\mathcal{B}_\SSS(r) = \{0\}$ for all $g \in M$;
\item[(iii)] There exists $r > 0$ such that $g\OO^n \cap B_\SSS(r) = \{0\}$ for all $g \in
M$ .
\end{enumerate}
\end{thm}

The equivalence between (i) and (iii) is proved in \cite[Theorem
5.12]{KlTo} and it is obvious that (ii) implies (iii). In order to
prove that (iii) implies (ii) note that, in view of the formula
(\ref{Artin}), every $\mathcal{B}_\SSS(r)$ is invariant under the
multiplication by elements from $\OO^*$. Now the implication easily
follows from the following lemma:

\begin{lem}
\label{balans} There exists a constant $\kappa
> 1$ with the following property. Let  $\textbf{x} =
(\textbf{x}^{(v)})_{v \in \SSS} \in K_\SSS^n$ be such that
$\textbf{x}^{(v)} \neq 0$ for all $v \in \SSS$. For each $v \in
\SSS$ we choose a positive real number $a_v$ in such a way that
$\mathbf{c}_\SSS(\textbf{x}) = \prod_{v \in \SSS}a_v$. Then there
exists $\xi \in
 \OO^*$ such that

\begin{equation}
\label{reduction1}\frac{a_v}{\kappa} \leq \|\xi \textbf{x}^{(v)}\|_v
\leq \kappa a_v\;
\end{equation}
for all $v \in \mathcal{S}$. In particular, for every $\textbf{x}$
as above there exists $\xi \in
 \OO^*$ such that
\begin{equation}
\label{reduction2} \frac{\mathbf{c}_\SSS(\textbf{x})^{1/m}}{\kappa}
\leq \|\xi \textbf{x}\|_\SSS \leq
\kappa\mathbf{c}_\SSS(\textbf{x})^{1/m},
\end{equation}
where $m = \# \SSS$.
\end{lem}

{\bf Proof.} Let $K_\SSS^1 = \{y = (y^{(v)}) \in K_\SSS^*| \prod_{v
\in \SSS}|y^{(v)}|_v = 1\}$. Then $\OO^* \subset K_\SSS^1$ and
$K_\SSS^1/\OO^* $ is compact  \cite[ch.2, Theorem 16.1]{CF}.
Therefore there exists a constant $\kappa_0 > 1$ such that for every
$y = (y^{(v)}) \in K_\SSS^1$ there exists $\xi \in
 \OO^*$ such that

\begin{equation}
\label{reduction3} \frac{1}{\kappa_0} \leq |\xi{y}^{(v)}|_v \leq
\kappa_0, \forall v \in \mathcal{S}.
\end{equation}

Let $\textbf{x}$ and $a_v, v \in \SSS$, be as in the formulation of
the proposition. There exists a constant $c > 1$, depending only on
$\SSS$, such that for every $v \in \SSS$ there exists $\alpha^{(v)}
\in K_v^*$ with
\begin{equation}
\label{reduction3'} \frac{c}{| \alpha^{(v)} |_v} \leq a_v \leq c|
\alpha^{(v)} |_v
\end{equation}
and $\prod_{v \in \SSS}| \alpha^{(v)} |_v = \prod_{v \in \SSS}a_v$.
So, $\mathbf{c}_\SSS(\alpha^{-1}\mathbf{x}) = 1$ where $\alpha =
(\alpha^{(v)})_{v \in \SSS} \in K_\SSS^*$. Put $\kappa = \kappa_0
c$. In view of (\ref{reduction3}) and (\ref{reduction3'}) there
exists $\xi \in \OO^*$ such that
$$
\frac{|\alpha^{(v)}|_v}{\kappa} \leq |\xi\mathbf{x}^{(v)}|_v \leq
\kappa |\alpha^{(v)}|_v, \forall v \in \mathcal{S},
$$
which proves (\ref{reduction1}).

In order to prove (\ref{reduction2}) it is enough to apply
 (\ref{reduction1}) with $a_v = \mathbf{c}_\SSS(\textbf{x})^{1/n}$.
\qed

\subsection{Horospherical subsets}
We need to prove a compactness criterion which reflects the group
structure of $\G$. %To this end we will
%generalize some propositions from our paper \cite{Tomanov-Weiss}.

We generalize the notion of horospherical subset from
\cite[Definition 3.4]{Tomanov-Weiss}.

\begin{defn}
\label{S-horospherical} Let $\RR \subset \SSS$. A finite subset
$\mathcal{M}$ of $\mathfrak{g}_\RR$ is called {\it
$\RR$-horospherical} (or, simply, horospherical when $\RR$ is
implicit) if $\mathcal{M} = \mathrm{pr}_\RR(\Ad(g)(\mathcal{M}_0))$,
where $g \in G$ and $\mathcal{M}_0$ is a subset of $\LLL$ which
spans linearly the unipotent radical of a maximal parabolic
$K$-subalgebra of $\LL$.
\end{defn}

The next proposition provides a compactness criterion in terms of
the intersection of pseudo-balls (and balls) in $\mathfrak{g}$ with %$\OO$-modules
$\mathfrak{g}_x$, $x \in
 G/\Gamma$ (see 2.1 for the notation). It generalizes \cite[Propositions 3.3 and
 3.5]{Tomanov-Weiss}.

\medskip
\begin{prop}
\label{intersection1} Assume that $\G$ is a semisimple algebraic
group. Then the following assertions hold:
\begin{enumerate}
 \item[(a)] There exists $r > 0$  $($respectively, $t > 0$$)$ such that for any $x = \pi(g)$ the
subalgebra of $\mathfrak{g}$ spanned by
%$\mathrm{Ad}_\SSS(g^{-1})(\mathcal{B}(r_0) \cap \mathcal{G}_x)$
$\mathcal{B}_\SSS(r) \cap \mathfrak{g}_x$ $($respectively,
${B}_\SSS(t) \cap \mathfrak{g}_x$$)$ is unipotent;

 \item[(b)]$\mathbf{(Compactness \ Criterion)}$ A subset $M$ of $G/\Gamma$ is relatively compact
 if and only if there exists $r > 0$ $($respectively, $t > 0$$)$ such that $\mathcal{B}_\SSS(r) \cap \mathfrak{g}_x$
$($respectively, ${B}_\SSS(t) \cap \mathfrak{g}_x$$)$ does not
contain a horospherical subset for any $x \in M$.
\end{enumerate}

\end{prop}

\subsection{Proof of Proposition \ref{intersection1}}
%\subsubsection{}
%\label{balls}
For every $t > 0$ we let $r = \big(\frac{t}{\kappa}\big)^m$, where
$\kappa$ and $m$ are as in the formulation of Lemma \ref{balans}. It
follows from Lemma \ref{balans} that
$$
{B}_\SSS({t}/{\kappa}) \subset \mathcal{B}_\SSS(r) \subset
\mathcal{O}^*{B}_\SSS(t).
$$
Now the validity of the proposition for the balls $B_\SSS(t)$
implies easily its validity for the pseudoballs
$\mathcal{B}_\SSS(r)$.

Further on, the proof of the proposition breaks in two cases. (In
view of 2.4, we will assume that $\Gamma = \G(\OO)$.)

\subsubsection{The case $\SSS = \SSS_\infty$}
\label{infty} Let $\rm{R}_{K/\Q}$ be the Weil restriction of scalars
functor. Then $\mathbf{H} = \rm{R}_{K/\Q}(\G)$ is a semisimple
$\Q$-algebraic group and $\HH = \rm{R}_{K/\Q}(\LL)$ is its $\Q$-Lie
algebra. Denote $\Delta = \mathbf{H}(\Z)$, $H  = \mathbf{H}(\R)$ and
$\mathfrak{h} = \HH(\R)$. The following properties of the functor
$\rm{R}_{K/\Q}$ are well known and easily follow from its definition
(see, for example, \cite[ch.2, \S 2.1.1]{PlaRa}). There exist
continuous isomorphisms $\mu:G \rightarrow H$ and $\nu: \mathfrak{g}
\rightarrow \mathfrak{h}$  such that $\mu(\Gamma) = \Delta$,
$\nu(\LL(\OO)) = \HH(\Z)$ and
$$\nu(\mathrm{Ad}_{G}(g)x) =
\mathrm{Ad}_{H}(\mu(g))\nu(x)$$ for all $g \in G$ and $x \in
\mathfrak{g}$. Moreover, $\nu$ maps bijectively the family of the
horospherical subsets of $\mathfrak{g}$ to the family of the
horospherical subsets of $\mathfrak{h}$ and $\mu$ induces an
homeomorphism $G/\Gamma \rightarrow H/\Delta$. Hence, when $\SSS =
\SSS_\infty$ the proposition follows from the case $K = \Q$
considered in \cite[Propositons 3.3 and 3.5]{Tomanov-Weiss}.

\subsubsection{The case $\SSS \varsupsetneq \SSS_\infty$}
\label{general} We introduce the topological rings ${\OO}_f \bydefn
\prod_{v \in \mathcal{S}_f}{\OO}_v$ and $K_f \bydefn K_\infty \times
\OO_f$ (see 2.1). So, $\OO_\infty = \OO \cap (K_\infty \times
\OO_f)$ is the ring of integers of $K$.

If $\widetilde{\G}$ is the simply connected covering of the
algebraic group $\G$ then $\widetilde{G}/\widetilde{\Gamma}$ is
naturally homeomorphic to $G/\Gamma$, where $\widetilde{G} =
\widetilde{\G}(K_\SSS)$ and $\widetilde{\Gamma} =
\widetilde{\G}(\OO)$. In view of this and of Theorem \ref{volumes}
below, we may (and will) assume without loss of generality that $\G$
is simply connected and without $K$-anisotropic factors. Then the
diagonal imbedding of $\Gamma$ into $\prod_{v \in \SSS_f}\G(K_v)$ is
dense. (This fact follows immediately from the strong approximation
theorem \cite[Theorem 7.12]{PlaRa}.) Therefore
$$
G = \G(K_f)\Gamma.
$$
 %Since $\G(K_f) = G_\infty
%\times \G(\OO_f)$,
Every $g \in G$ can be writhen in the following way
\begin{equation}
\label{decomp} g = g_\infty g_f \gamma,
\end{equation}
where $g_\infty \in G_\infty, g_f \in \G(\OO_f)$ and $\gamma \in
\Gamma$. Let $\Gamma_\infty = \G(K_f) \cap \Gamma$. Then $G/\Gamma$
is homeomorphic to $\G(K_f)/\Gamma_\infty$ and the projection of
$\G(K_f)$ on $G_\infty$  yields the following map
 $$\varphi:G/\Gamma \rightarrow G_\infty/\Gamma_\infty, \
\varphi(\pi(g)) \bydefn \pi_\infty(g_\infty), \forall g \in G,$$
where $\pi_\infty: G_\infty \rightarrow G_\infty/\Gamma$ is the
natural map. In view of the compactness of $\G(\OO_f)$, $\varphi$ is
a proper continuous map.

Let $A$ be a subset of $\mathfrak{g}_\infty$ and $x = \pi(g)$ for
some $g \in G$. Set $A_f = A \times \LL(\OO_f)$. Using
(\ref{decomp}) and the fact that $\LL(\OO)$ is invariant under the
adjoint action of $\Gamma$, we obtain
\begin{equation}
\begin{split}
\label{proj} \mathrm{pr}_\infty\big(\mathfrak{g}_x \cap A_f\big) =
\mathrm{pr}_\infty\Big(\Ad_{\SSS}&(g)\big(\LL(\OO)\big) \cap A_f\Big) =  \\
%\Ad_{\SSS}(g_\infty g_f)\big(\LL(\OO)\big) \cap \LL(K_f)& =
\mathrm{pr}_\infty\Big(\Ad_{\SSS}(g_\infty g_f)\big(\LL(\OO) \cap
&A_f\big)\Big)  =
%\Ad_{\infty}(g_\infty)\big(\LL(\OO_\infty)\big) \cap A =
\mathfrak{g}_{\infty,y} \cap A,
\end{split}
\end{equation}
where $y = \varphi(x)$ and $\mathfrak{g}_{\infty,y} = \Ad_\infty
(g_\infty) \LL(\OO_\infty)$. (Recall that $\mathrm{pr}_\infty$
denotes the natural projection $\mathfrak{g} \rightarrow
\mathfrak{g}_\infty$.)

Let $\tilde{B}(t) = B_\infty(t) \times
\LL(\OO_f)$. Applying (\ref{proj}) with $A = B_\infty(t)$, we get
%for any $x \in G/\Gamma$
$$
\mathrm{pr}_\infty\big(\mathfrak{g}_x \cap \tilde{B}(t)\big) =
\mathfrak{g}_{\infty, y} \cap B_\infty(t).
$$
Since the restriction of $\mathrm{pr}_\infty$ to $\mathfrak{g}_x$ is
injective, we obtain that the subalgebra spanned by $\mathfrak{g}_x
\cap \tilde{B}(t)$ is unipotent if and only if the subalgebra
spanned by $\mathfrak{g}_{\infty, y} \cap B_\infty(t)$ is unipotent.
This, in view of 3.3.1, proves (a).

Let us prove (b). If $M$ is compact, it follows from the continuity
of the adjoint action that if $t > 0$ is sufficiently small then
$B_\SSS(t) \cap \mathfrak{g}_x$ does not contain horospherical
subsets for all $x \in G/\Gamma$. In order to prove the inverse
implication, let $M \subset G/\Gamma$ and $t > 0$ be such that
$B_\SSS(t) \cap \mathfrak{g}_x$ does not contain horospherical
subsets for any $x \in M$. Assume the contrary, that is, that there
exists  a divergent sequence $\{x_i\}$ of elements in $M$. Then the
sequence $\{y_i = \varphi(x_i)\}$ is divergent in
$G_\infty/\Gamma_\infty$ (because $\varphi$ is proper). Since the
proposition is true for $G_\infty/\Gamma_\infty$, for every
$\varepsilon > 0$ there exists $i \gg 0$ such that
$B_\infty(\varepsilon) \cap \mathfrak{g}_{\infty,y_i}$ contains a
horospherical subset. Set $\tilde{B}(\varepsilon) =
B_\infty(\varepsilon) \times \LL(\OO_f)$. By (\ref{proj}) (applied
with $A = B_\infty(\varepsilon)$) and the injectivity of the
restriction of $\mathrm{pr}_\infty$  to $\mathfrak{g}_x$, we obtain
that $\tilde{B}(\varepsilon) \cap \mathfrak{g}_{x_i}$ contains a
horospherical subset. Now, using Lemma \ref{balans}, we conclude
that $B_\SSS(t) \cap \mathfrak{g}_{x_i}$ contains horospherical
subsets for all sufficiently large $i$. Contradiction. \qed

\subsection{Expanding transformations} For every $v \in \SSS$, we fix a maximal
$K_v$-split torus $\mathbf{T}_v$ of $\G$. We denote $T_v =
\T_v(K_v)$ and $T_\mathcal{R} = \prod_{v \in \RR}T_v$ where
$\mathcal{R}$ is a non empty subset of $\SSS$.

\begin{prop} \label{pushing out} With the above notation,  for every real $\tau > 1$
there exists a finite set $t_1, \ldots, t_{s}$ of elements in
$T_\RR$ such that if $\mathfrak{u}$ is a unipotent subalgebra of
$\mathfrak{g}_\RR$ then there exists an element $t_i$ such that

\begin{equation}
\label{expansion1} \|\Ad(t_i) (\mathbf{x})\|_\RR \geq \tau
\|\mathbf{x}\|_\RR
\end{equation}
for all $\mathbf{x} \in \mathfrak{u}$.
\end{prop}

 {\bf Proof.} It is easy to see that it is enough to prove the proposition when $\RR$ is a
singleton. Let $\RR = \{v \}$. If $v$ is real then the proposition
is proved in \cite[Proposition 4.1]{Tomanov-Weiss}. Here we present
a shorter proof for an arbitrary $v$.

Let  $\mathfrak{u}^+_v$ and $\mathfrak{u}^-_v$ be invariant under
the adjoint action of $T_v$ maximal unipotent subalgebras  of
$\mathfrak{g}_v$ which are opposite to each other.
%(see\cite{Borel}).
Then
\begin{equation}
\label{chamber} C_v  =  \{d \in T_v| \lim_{n \rightarrow
+\infty}\Ad(d^n)x = \infty, \forall x \in \mathfrak{u}^+_v\}\\
\end{equation}
is the interior of the Weil chamber corresponding to
$\mathfrak{u}^+_v$ (see \cite{Borel}). Denote by $U_v^+$ and $U_v^-$
the unipotent subgroups of $G_v$ with Lie algebras
$\mathfrak{u}^+_v$ and $\mathfrak{u}^-_v$, respectively.

Now let $\mathfrak{u}_v$ be any maximal unipotent subalgebra of
$\mathfrak{g}_v$. There exists $g \in G_v$ such that
$\Ad(g)\mathfrak{u}^+_v = \mathfrak{u}_v$. By Bruhat decomposition
$g = a\omega b$, where $\omega \in \mathcal{N}_{G_v}(T_v)$, $a$ and
$b \in U_v^+$ and $\omega^{-1}a\omega \in U_v^-$. We can write
$\mathfrak{u}_v = \Ad(\omega a^-)\mathfrak{u}^+_v$, where $a^- =
\omega^{-1}a\omega$. Let $x \in \mathfrak{u}_v^+$ and $ f_v \in
C_v$. We put $y = \Ad(\omega a^-)x$ and  $d_v =
 \omega f_v \omega^{-1}$. Using (\ref{chamber}) and the fact that $\lim_{n
\rightarrow +\infty} f_{v}^n a^-f_{v}^{-n} = 0$,  we get
\begin{eqnarray*}
\lim_{n \rightarrow +\infty}\Ad(d_{v}^n)y = \lim_{n \rightarrow
+\infty}\Ad(\omega(f_{v}^n a^-f_{v}^{-n}))\circ \Ad(f_{v}^n)(x) =
\infty.
\end{eqnarray*}
Therefore, taking $t = d^n$ with $n$ sufficiently large, we obtain
that
\begin{equation}
\label{expansion2}
 \Vert \Ad(t)z \Vert_v > \tau \Vert z \Vert_v
\end{equation}
for all non-zero $z \in \mathfrak{u}_v$.

Since the stabilizer of every maximal unipotent subalgebra is a
minimal parabolic subgroup and all minimal parabolic subgroups are
conjugated, the set of all maximal unipotent subalgebras can be
identified with the compact homogeneous space $G_v/P^{+}_v$, where
$P^{+}_v$ is the parabolic subgroup of $G_v$ with Lie algebra
$\mathfrak{u}^{+}_v$. It is easy to see that (\ref{expansion2}) is
true for all subalgebras in a neighborhood of $\mathfrak{u}_v$. Now
the existence of the elements $t_1, \ldots, t_{s}$ as in the
formulation of the theorem follows from the compactness of
$G_v/P^{+}_v$ by a standard argument. \qed

\section{Closed orbits of reductive $K$-groups}
\medskip

\subsection{Reductive groups}
Recall the $\SSS$-adic version of a well-known theorem of Borel and
Harish-Chandra. (As usual, $G = \G(K_\SSS)$ and $\Gamma = \G(\OO)$.)

\medskip
\begin{thm}
\label{volumes} \mbox{\rm{(}}cf.\cite[Theorem 5.7
 ]{PlaRa}\mbox{\rm{)}} Let $\G$ be a reductive $K$-group and let
 $\mathrm{X}_K(\G)$ be the group of $K$-rational characters of $\G$.
 Then
 \begin{enumerate}
\item[(a)] $G/\Gamma$ has a finite invariant volume if and only if $\mathrm{X}_K(\G) =
\{1\}$;
\item[(b)] $G/\Gamma$ is compact if and only if $\G$ is anisotropic over
$K$.
\end{enumerate}
\end{thm}
\medskip
Because of the lack of appropriate reference we will prove the
following known proposition.

\medskip
\begin{prop}
\label{closed orbits1} With the above notation, let $\mathbf{H}$ be
a reductive subgroup of $\G$ defined over $K$ and $H =
\mathbf{H}(K_\SSS)$. Then $H\pi(e)$ is closed in $G/\Gamma$.
\end{prop}
\medskip

{\bf Proof.} Using the Weil restriction of scalars, one can reduce
the proof to the case when $K = \Q$. In view of \cite[Proposition
7.7]{Borel - arithmetique} there exists a $\Q$-rational action of
$\G$ on an affine $\Q$-variety $\mathbf{V}$ admitting an element $a
\in \mathbf{V}(\Z)$ such that $\mathbf{H} = \{g \in \G|ga = a\}$.
Since the map $\G \rightarrow \mathbf{V}$, $g \rightarrow ga$, is
polynomial with rational coefficients, there exists a non-zero
integer $n$ such that $\gamma n a \in \mathbf{V}(\OO)$ for all
$\gamma \in \Gamma$. Therefore $\Gamma H$ is closed in $G$,
equivalently, $H \pi(e)$ is closed. \qed

\medskip
\subsection{Algebraic tori}
We will need the following
\begin{prop}
\label{K-torus} Let $\mathbf{T}$ be a $K$-torus in $\G$ and let
$\mathcal{R}$ be a non-empty subset of $\SSS$. Suppose that $T_\RR$
is not compact. Then the orbit $T_\RR \pi(e)$ is divergent if and
only if the following conditions are fulfilled:
\begin{enumerate}
\item[(i)] $\RR = \{v_\circ\}$ is a singleton, and
\item[(ii)] $\mathrm{rank}_K \mathbf{T} = \mathrm{rank}_{K_{v_\circ}} \mathbf{T} > 0$.
\end{enumerate}
\end{prop}
{\bf Proof.} In view of Proposition \ref{closed orbits1}  the orbit
$\mathbf{T}(K_\SSS) \pi(e)$ is closed and, therefore, homeomorphic
to $\mathbf{T}(K_\SSS)/(\mathbf{T}(K_\SSS) \cap \Gamma)$. So, we may
suppose, with no loss or generality, that $\mathbf{T} = \G$.

Assume that the orbit $T_\RR \pi(e)$ is divergent. Let $\te_a$
(respectively, $\te_d$) be the largest $K$-anisotropic
(respectively, split over $K$) subtorus of $\te$. It is well known
that $\te$ is an almost direct product of $\te_a$ and $\te_d$. This
implies that if there exists $v \in \RR$ such that
$\mathrm{rank}_{K_v} \mathbf{T} > \mathrm{rank}_{K} \mathbf{T}$ then
$\mathbf{T}_a(K_\RR)$ is not compact. But
$\mathbf{T}_a(K_\SSS)\pi(e)$ is compact (Theorem \ref{volumes}).
Therefore, $T_\RR \pi(e)$ can not be divergent, a contradiction. So,
$\mathrm{rank}_{K_v} \mathbf{T} = \mathrm{rank}_{K} \mathbf{T}$ for
all $v \in \RR$. In this case $\mathbf{T}_a(K_\RR)$ is compact and,
since  $T_\RR$ is not compact, $\te_d$ is not trivial. Note that
$T_\RR \pi(e)$ is divergent if and only if
$\mathbf{T}_d(K_\RR)\pi(e)$ is divergent.

In order to prove (i) consider the character group
$\mathrm{X}_K(\mathbf{T})$ of $\te$.  It is well known that
$\mathrm{X}_K(\mathbf{T})$ is a free $\Z$-module of rank equal to
$\dim \te_d$ (cf. \cite[8.15]{Borel}). Let $\chi_1, \ldots, \chi_r$
be a basis of $\mathrm{X}_K(\mathbf{T})$. Define a homomorphism of
$K$-algebraic groups $\chi = (\chi_1, \ldots, \chi_r): \mathbf{T}
\rightarrow \mathbf{G}_m^r$, where $\mathbf{G}_m$ denotes the
one-dimensional $K$-split torus. Let $T = \T(K_\SSS)$  and $T_\circ
= \{(t_v)_{v \in \SSS} \in {T} |\prod_{v \in \SSS}|\chi_i(t_v)|_v =
1 \ \mathrm{for \ all \ }i \}$. It follows from \cite[ch.2, Theorem
16.1]{CF} that $\Gamma$ is a co-compact lattice in $T_\circ$. Set
$\varphi: T \rightarrow \R^r$, $\varphi((t_v)_{v \in \SSS}) =
\big(\log(\prod_{v \in \SSS}|\chi_1(t_v)|_v), \ldots,$
$\log(\prod_{v \in \SSS}|\chi_r(t_v)|_v) \big)$. It is clear that
$\varphi$ is a continuous surjective homomorphism of locally compact
topological groups with $\ker(\varphi) = T_\circ$. Since
$T_\circ/\Gamma$ is compact, $\varphi$ induces a proper homomorphism
$\psi: T/\Gamma \rightarrow T/T_\circ$. Now let $\RR$ contain two
different valuations $v_1$ and $v_2$. It is easy to find sequences
$\{a_i\}$ in $K_{v_1}^*$ and $\{b_i\}$ in $K_{v_2}^*$ such that
$\log |a_i|_{v_1} \rightarrow +\infty$, $\log |b_i|_{v_2}
\rightarrow -\infty$ and the sequence $\{\log |a_i|_{v_1} + \log
|b_i|_{v_2} \}$ is bounded. We define a sequence $\{s_i =
(s^{(v)}_i)_{v \in \RR}\}$ in $T_\RR$ as follows:

\begin{displaymath}
s_i^{(v)} = \left\{ \begin{array}{ll}
1, \ \textrm{if} \  v \in \RR \setminus \{v_1, v_2\};\\
\chi_1(s_i^{(v_1)}) = a_i \ \textrm{and} \ \chi_j(s_i^{(v_1)}) = 1 \
\textrm{for all} \ j > 1;\\
\chi_1(s_i^{(v_2)}) = b_i \ \textrm{and} \ \chi_j(s_i^{(v_2)}) = 1 \
\textrm{for all} \ j > 1.
\end{array} \right.
\end{displaymath}
We have that $\{s_i\}$ is unbounded and that $\{\varphi(s_i)\}$ is
bounded. (Recall that $T_\RR$ is considered as a subgroup of $T$, so
that the notation $\varphi(s_i)$ makes sense.) Since $\psi$ is
proper, $s_i \pi(e)$ is bounded. Therefore the orbit $T_\RR \pi(e)$
is not divergent. This contradiction completes the proof of (i).

Assume that $\RR$ contains only one valuation $v_\circ$ and that
$\mathrm{rank}_K \mathbf{T} = \mathrm{rank}_{K_{v_\circ}} \mathbf{T}
> 0$. It follows from the above definition of $\varphi$ and the
fact that $\chi$ is an homomorphism with compact kernel, that if a
sequence $\{t_i\}$ in $T_\RR$ diverges then $\{\varphi(t_i)\}$ does
too. Therefore $T_\RR \pi(e)$ is a divergent orbit. \qed

\medskip

Proposition \ref{K-torus} implies:

\begin{prop}
\label{K-torus1} Let $\mathbf{T}$ be a $K$-torus and let
$\mathcal{R}$ be a non-empty subset of $\SSS$. Then the orbit $T_\RR
\pi(e)$ is closed if and only if one of the following conditions
holds:
\begin{enumerate}
\item $\RR = \SSS$;
\item $\mathrm{rank}_{K_{v}} \mathbf{T} = 0$ for all $v \in \RR$,
equivalently, $T_\RR$ is compact;
\item $\RR = \{v_\circ\}$ and $\mathrm{rank}_K \mathbf{T} = \mathrm{rank}_{K_{v_\circ}} \mathbf{T}$.
\end{enumerate}
\end{prop}
{\bf Proof.} Note that if $\RR \neq \SSS$ and $T_\RR$ is not compact
then $T_\RR \pi(e)$ is closed if and only if it is divergent. Now
the proposition follows easily from Proposition \ref{K-torus}.
\qed

\medskip

\section{Parabolic subgroups and divergent orbits}
\medskip

\subsection{Main proposition}  Recall that, given a subset $\mathcal{R} \subset
\mathcal{S}$, we use the notation $\mathrm{pr}_\mathcal{R}$ to
denote depending on the context the projection $G \rightarrow
G_\mathcal{R}$ or the projection $\mathfrak{g} \rightarrow
\mathfrak{g}_\mathcal{R}$.

The goal of this section is to prove the following

\medskip

\begin{prop}
\label{split torus} Let $\G$ be a reductive  $K$-algebraic group,
$\mathcal{R}$ be a non-empty subset of $ \SSS$,  $g =(g_v)_{v \in
\SSS} \in G$ and $x = \pi(g)$. Assume that $\rank_K \G > 0$ and that
for every minimal parabolic $K$-subalgebra $\bb$ of $\LL$ containing
the Lie algebra of $\mathbf{D}$ there exists a horospherical subset
$\mathcal{M}_\mathfrak{b}$ of $\mathfrak{g}_\mathcal{R}$ such that
$\mathcal{M}_\mathfrak{b} \subset
\mathrm{pr}_\mathcal{R}(\mathfrak{g}_x) \cap
\mathfrak{b}_{\mathcal{R}}$. Then the following assertions hold:
\begin{enumerate}
\item[(a)] For every $v \in \RR$ the orbit $D_v \pi(g)$ is divergent;
\item[(b)] If $g_{\RR} = \mathrm{pr}_\RR(g)$ then
\begin{equation}
 \label{int0'} g_{\RR} \in
\mathcal{Z}_{G_{\RR}}(D_{\RR}) \mathrm{pr}_\mathcal{R}(\G(K));
\end{equation}
%where $\G(K)$ is identified with $\mathrm{pr}_\mathcal{R}(\G(K))$;
\item[(c)] There exists a maximal $K$-split torus
%$\mathbf{D}^\bullet$
$\mathbf{S}$ of $\G$ such that
%$$D^\bullet_v = {g_v}^{-1}D_vg_v$$
\begin{equation}
 \label{int10}
S_v = {g_v}^{-1}D_vg_v
\end{equation}
for all $v \in \mathcal{R}$, where $S_v = \se(K_v)$.
\end{enumerate}
\end{prop}

\medskip

In order to prove Proposition \ref{split torus} we will need some
facts from algebraic group theory.

\subsection{Intersections of parabolic subgroups} The next three propositions remain valid for any field
$K$.

\begin{prop}
\label{prop6}
 \cite[Propositions 14.22 and 21.13]{Borel} Let $\pe$ and  $\mathbf{Q}$ be  parabolic
$K$-subgroups of $\G$.
\begin{enumerate}
\item[(i)] $(\mathbf{P} \cap \mathbf{Q}){R}_u(\pe)$ is a parabolic
$K$-subgroup;
\item[(ii)]  If $\mathbf{Q}$ is conjugate to $\pe$ and contains
${R}_u(\pe)$ then $\mathbf{Q}= \pe$.
\end{enumerate}
\end{prop}

We also have

\begin{prop}
\label{prop77}\cite[Proposition 5.2]{Tomanov-Weiss} For every
minimal parabolic $K$-subgroup $\mathbf{B}$ containing $\mathbf{D}$
we let $\mathbf{P}_{\mathbf{B}}$ be a proper parabolic $K$-subgroup
containing $\mathbf{B}$. Then

\begin{equation}
\label{eq: prop7} \bigcap_{\mathbf{B}} \mathbf{P_B} =
\mathcal{Z}_\G(\mathbf{D}).
\end{equation}
\end{prop}

\medskip
Keeping the notation and assumptions of Proposition \ref{prop77}, we
prove:

\begin{prop}
\label{prop7'} Let $n \in \mathcal{N}_\G(\mathcal{Z}_\G
(\mathbf{D}))$. Assume that for every $\mathbf{B}$ the group
$n\mathbf{P_B}n^{-1}$ is defined over $K$. Then $n \in
\mathcal{N}_\G (\mathbf{D})$. The projection of $n$ into the Weyl
group $\mathrm{W}_K = \mathcal{N}_\G (\mathbf{D})/\mathcal{Z}_\G
(\mathbf{D})$ is uniquely defined by  the map $\mathbf{B}
\rightarrow n\mathbf{P_B}n^{-1}$.

\end{prop}

{\bf Proof.} The uniqueness of the projection of $n$ into
$\mathrm{W}_K$ follows immediately from Proposition \ref{prop77} and
the fact that every parabolic subgroup coincides with its
normalizer.

We will assume that for every $\be$ the group $\pe_{\be}$ is minimal
among the parabolic $K$-subgroups $\pe$ containing $\be$ and such
that $n \pe n^{-1}$ is defined over $K$.

Assume that there exists $\be$ such that $\be = \pe_{\be}$. Let
$\be' = n \be n^{-1}$. Since all minimal parabolic $K$-subgroups are
conjugated under the action of $\mathrm{W}_K$ and
$\mathcal{N}_\G(\mathbf{D}) =
\mathcal{N}_\G(\mathbf{D})(K)\mathcal{Z}_\G (\mathbf{D})$
\cite[Theorem 21.2]{Borel}, there exists $n_\circ \in
\mathcal{N}_\G(\mathbf{D})(K)$ such that $\be = n_\circ
\be'n_\circ^{-1}$. Therefore, $\be = n_\circ n \be (n_\circ n)^{-1}$
which implies that $n_\circ n \in \be$. Since
$\mathcal{N}_\G(\mathbf{D}) \subset \mathcal{N}_\G(\mathcal{Z}_\G
(\mathbf{D}))$, we get $n_\circ n \in \mathcal{N}_\be(\mathcal{Z}_\G
(\mathbf{D}))$. Now, the proposition follows from the fact that
$\mathcal{Z}_\G (\mathbf{D}) = \mathcal{N}_\be(\mathcal{Z}_\G
(\mathbf{D}))$ \cite[Corollary 14.19]{Borel}.

Assume that $\pe_\be \varsupsetneqq \be$ for all $\be$. Choose a
$\pe_\be$ with the minimal dimension and set $\pe = \pe_\be$. Let
$\Phi(\de, \G)$ be the relative root system of $\G$ with respect to
$\de$. (See \cite[21.1 and 8.17]{Borel} for the standard definition
of a system of $K$-roots.) Since $\pe \varsupsetneqq \be$, there
exists  a long root $\alpha \in \Phi(\de, \G)$ such that $\pm\alpha$
are roots of the group $\pe$ with respect to $\de$. Recall that all
roots of the same length in $\Phi(\de, \G)$ are conjugated under the
action of $\mathrm{W}_K$ \cite[10.4, Lemma C and 10.3,
Theorem]{Humphreys-Lie}. Therefore there exists a minimal parabolic
$K$-subgroup $\be^+$ containing $\de$ such that $\alpha$ is a
maximal long root of $\be^+$ relative to $\de$. Let $\Delta^+$ be
the set of simple roots corresponding to $\be^+$. Then in the
expression of $\alpha$ as a linear combination of the roots in
$\Delta^+$ all coefficients are strictly positive \cite[10.4, Lemma
A]{Humphreys-Lie}. It follows from the explicit description of the
standard parabolic $K$-subgroups (see \cite[21.11]{Borel}), that
$-\alpha$ is not a root of any parabolic $K$-subgroup containing
$\be^+$. Similarly, $\alpha$ is not a root of any parabolic
$K$-subgroup containing $\be^-$, where $\be^-$ is the minimal
parabolic $K$-subgroup opposite to $\be^+$. As a consequence, one of
the $K$-subgroups $(\pe_{\be^+} \cap \mathbf{P}){R}_u(\pe)$ or
$(\pe_{\be^-} \cap \mathbf{P}){R}_u(\pe)$ is strictly smaller than
$\mathbf{P}$. Let $\pe \neq (\pe_{\be^+} \cap
\mathbf{P}){R}_u(\pe)$. Since $(\pe_{\be^+} \cap
\mathbf{P}){R}_u(\pe)$ is a parabolic $K$-subgroup (Proposition
\ref{prop6}(i)) and $n (\pe_{\be^+} \cap
\mathbf{P}){R}_u(\pe)n^{-1}$ is defined over $K$. The latter
contradicts the choice of $\pe$, which completes our proof. \qed

\medskip

\begin{remark}
\label{quasisplit group} In connection with the above proposition,
let us note that in certain cases $\mathcal{N}_\G(\mathbf{D})
\varsubsetneq \mathcal{N}_\G(\mathcal{Z}_\G (\mathbf{D}))$. As a
simple example one can consider the special unitary group
$\mathbf{SU}_3(h)$, where $h$ is an hermitian form with coefficients
from $K$ of indice 1. This is a quasisplit group of type $A_2$.
Therefore $\mathcal{N}_\G(\mathcal{Z}_\G
(\mathbf{D}))/\mathcal{Z}_\G (\mathbf{D})$ is isomorphic to the
symmetric group $S_3$ and $\mathcal{N}_\G(\mathbf{D})/\mathcal{Z}_\G
(\mathbf{D})$ is a group of order two.
\end{remark}

\medskip

\subsection{Proof of Proposition \ref{split torus}}
We start the proof with a general remark. We keep the notation from
the formulation of the proposition. For every $\bb$ there exists a
finite subset $\mathcal{M}_\mathfrak{b}^\bullet$ of $\LLL$ which
spans linearly the unipotent radical of a maximal parabolic
$K$-subgroup $\mathbf{P}_\mathfrak{b}^\bullet$ of $\G$ and such that
$\mathcal{M}_\mathfrak{b} =
\mathrm{pr}_\RR(\Ad(g)(\mathcal{M}_\mathfrak{b}^\bullet)$. So, if $v
\in \mathcal{R}$, we have
$$
g_v {R}_u(\mathbf{P}_\mathfrak{b}^\bullet)(K_v)g_v^{-1} \subset
\mathbf{B}(K_v),
$$
where $\mathbf{B}$ is the $K$-algebraic subgroup of $\G$ the Lie
algebra of which is $\bb$. It follows from Proposition
\ref{prop6}(ii) that there exists a parabolic $K$-subgroup
$\mathbf{P}_\mathfrak{b}$ containing $\mathbf{B}$ such that
\begin{equation}
 \label{int3}
 \mathbf{P}_\mathfrak{b} = g_v \mathbf{P}_\mathfrak{b}^\bullet
 g_v^{-1}
\end{equation}
for all $v \in \mathcal{R}$.

Let us prove (a). (Remark that (a) follows {\it a posteriori} from
(b) and Proposition \ref{K-torus}.) Fix $v \in \mathcal{R}$. We want
to prove that the orbit $D_v\pi(g)$ diverges. Let $\{d_i\}$ be a
divergent sequence in $D_v$. Put $s_i = g_v^{-1}d_ig_v$. It is
enough to prove that the sequence $\{s_i\pi(e)\}$ is divergent.
Passing to a subsequence we may assume that $\{d_i^{-1}\}$ belongs
to the Weyl chamber corresponding to some minimal parabolic
$K$-subgroup $\mathbf{B}$. Let $\UU$ be the Lie algebra of
${R}_u(\mathbf{P}_{\mathfrak{b}}^\bullet)$. Let $m$ be the dimension
of $\UU$ and let $\bigwedge^{m}\Ad$ be the adjoint representation of
$\G$ on the $m$-th exterior power $\bigwedge^{m}\LL$. Since $\UU$ is
defined over $K$, there exists a non-zero $K$-{\it rational} vector
$z \in \bigwedge^{m}\LL$ corresponding to $\UU$. It is known (see
the proof of Proposition \ref{prop7'}) that if $\alpha$ is a maximal
root of $\mathbf{B}$ with respect to $\mathbf{D}$ then $\alpha$ is a
root of every standard parabolic subgroup containing $\mathbf{B}$
and, given the choice of $\{d_i\}$, $\lim_{i \rightarrow \infty}
\alpha (d_i) = 0$. Since $\mathbf{P}_{\mathfrak{b}} = g_v
\mathbf{P}_{\mathfrak{b}}^\bullet g_v^{-1}$ and
$\mathbf{P}_{\mathfrak{b}}$ is a parabolic containing $\mathbf{B}$,
we obtain that
$$\lim_{i \rightarrow \infty}\|\bigwedge^{m}\Ad(d_i)g_vz\|_v = 0.$$
%Therefore
%$$ \lim_{i \rightarrow \infty}\|\Ad(s_i)z\|_v \rightarrow 0.$$
This implies
$$
\lim_{i \rightarrow \infty}\mathbf{c}_\SSS(\bigwedge^{m}\Ad(s_i)z) =
0.
$$
It follows from Theorem \ref{mahler} (ii) that $\{s_i\pi(e)\}$
diverges. This completes the proof of (a).

Note that (c) follows immediately from (b). So, it remains to prove
(b). Let $\mathbf{P}_\mathfrak{b}^\bullet$ be as above. Set
$\mathbf{H} = \bigcap_\mathfrak{b} \mathbf{P}_\mathfrak{b}^\bullet$.
Since $\mathbf{P}_\mathfrak{b}$ is a $K$-parabolic subgroup of $\G$
containing $\mathbf{B}$, in view of Proposition \ref{prop77}, we get
that
\begin{equation}
\label{int0} \mathbf{H} = \bigcap_\mathfrak{b}
g_v^{-1}\mathbf{P}_\mathfrak{b}g_v = g_v^{-1}
\big(\bigcap_\mathfrak{b} \mathbf{P}_\mathfrak{b}\big) g_v =
g_v^{-1}\mathcal{Z_{\mathbf{G}}(\mathbf{D})}g_v
\end{equation}
for all $v \in \mathcal{R}$.

Note that the groups  $\mathcal{Z_{\mathbf{G}}(\mathbf{D})}$ and
$\mathbf{H}$ are reductive and defined over $K$. Let $\mathbf{Z}$
(respectively, $\mathbf{Z}^\bullet$) be the Zariski connected
component of the center of $\mathcal{Z_{\mathbf{G}}(\mathbf{D})}$
(respectively,  $\mathbf{H}$). It follows from (\ref{int0}) that
\begin{equation}
\label{int1} \mathbf{Z}^\bullet = g_v^{-1}\mathbf{Z}g_v
\end{equation}
for all $v \in \mathcal{R}$. Since  $\mathbf{D}$ is a maximal
$K$-split torus of $\G$, we have that $\mathbf{D} = \mathbf{Z}_d$,
where $\mathbf{Z}_d$ is the largest $K$-split subtorus  of
$\mathbf{Z}$.

Denote by $\mathbf{Z}^{\bullet}_d$ the largest $K$-split subtorus of
$\mathbf{Z}^\bullet$ and assume that $\mathbf{Z}^{\bullet}_d$ is not
maximal in $\G$. Let $\mathbf{Z}^{\bullet}_a$ be the largest
$K$-anisotropic subtorus of $\mathbf{Z}^\bullet$. Fix $v \in
\mathcal{R}$. Since every $K$-torus is an almost direct product over
$K$ of its largest $K$-split and its largest $K$-anisotropic subtori
\cite[Proposition 8.15]{Borel}, it follows from (\ref{int1}) that
there exists an element $t \in \mathbf{Z}^{\bullet}_a(K_v) \cap
g{_v}^{-1}\mathbf{D}(K_v)g_v$ such that $\{t^n | n \in \N\}$ is a
divergent sequence. In view of (a), $\{g_v t^n g{_v}^{-1}\pi(g)\}$,
and therefore $\{ t^n \pi(e)\}$, are also divergent sequences. The
latter contradicts the fact that the orbit
$\mathbf{Z}^{\bullet}_a(K_\mathcal{R})\pi(e)$ is compact (see
Theorem \ref{volumes}). Therefore $\mathbf{Z}^{\bullet}_d$ is a
maximal $K$-split torus of $\G$.

Since the maximal $K$-split tori are conjugated under $\G(K)$
\cite[Theorem 20.9]{Borel}, there exists $q \in \G(K)$ such that
$\mathbf{Z}^{\bullet}_d = q^{-1}\mathbf{D}q$. Also,
$\mathcal{Z}_\G(\mathbf{Z}^{\bullet}_d) =
q^{-1}\mathcal{Z}_{\G}(\mathbf{D})q$,
$\mathcal{Z}_\G(\mathbf{Z}^{\bullet}_d) \supset \mathbf{H}$ and
$\dim \mathbf{H} = \dim \mathcal{Z}_{\mathbf{G}}(\mathbf{D})$.
Therefore,
$$
\mathbf{H} = q^{-1}\mathcal{Z_{\mathbf{G}}(\mathbf{D})}q.
 $$
In view of (\ref{int0}),  we have
$$
g_v q^{-1} \in
\mathcal{N}_{\mathbf{G}}(\mathcal{Z}_{\mathbf{G}}(\mathbf{D})),
\forall v \in \mathcal{R}.
$$
Given $v \in \mathcal{R}$, the group
$$
qg_v^{-1}\mathbf{P}_\mathfrak{b}(qg_v^{-1})^{-1} = q
\mathbf{P}_\mathfrak{b}^\bullet q^{-1}
$$
is defined over $K$ for every $\mathfrak{b}$. It follows from
Proposition \ref{prop7'} that there exists $n \in
\mathcal{N}_\G(\mathbf{D})(K)$ such that
\begin{equation}
\label{int2} nqg_{v}^{-1} \in \mathcal{Z_{\mathbf{G}}(\mathbf{D})},
\forall v \in \mathcal{R}.
\end{equation}
Since $n$ is the same for all $v \in \mathcal{R}$, (\ref{int2})
implies (\ref{int0'}), which completes the proof. \qed

\medskip

\section{Proofs of Theorem \ref{thm2} and of its corollaries}

\subsection{Proof of Theorem \ref{thm2}.} Let the conditions
(i)-(iii) in the formulation of the theorem hold. Since $\rank_{K_v}
\G \geq \rank_K \G$, it follows from (ii) that $\rank_{K_v} \G =
\rank_K \G$ for all $v \in \RR$. Therefore, $T_\RR/D_\RR$ is
compact. So, $T_\RR\pi(g)$ is closed and locally divergent if and
and only if $D_\RR\pi(g)$ has this property. In view of (iii),
$g^{-1}D_\RR g = \widetilde{D}_\RR$, where $\widetilde{D}_\RR =
\widetilde{\de}(K_\RR)$ and $\widetilde{\de}$ is a $K$-split torus.
Using (i) and Proposition \ref{K-torus}, it is easy to see that
$\widetilde{D}_\RR \pi(g)$, and therefore ${D}_\RR \pi(g)$, are
closed locally divergent orbits.

Let the orbit $T_\RR\pi(e)$ be closed and locally divergent. In view
of Theorem \ref{volumes}(b), $\rank_K \G > 0$. Moreover, since every
$\te_v$ is a product of a maximal $K_v$-split torus and a compact,
we can suppose without loss of generality that $\te_v$ is a maximal
$K_v$-split torus.

Denote by $\mathbf{S}$ the connected component of the Zariski
closure of $g^{-1}{T}_{\RR}g \cap \Gamma$ in $\G$. Suppose that
$\mathbf{S}$ is not trivial. Then $\RR = \SSS$. Set $S =
\mathbf{S}(K_\SSS)$. Since $S$ is not compact, $S\pi(e)$ is locally
divergent and $\mathbf{S}$ is $K_v$-split, $v \in \SSS$, it follows
from Proposition \ref{K-torus} that $\mathbf{S}$ is $K$-split. Set
$\mathbf{H} = \mathcal{Z}_\G(\mathbf{S})$, $H = \mathbf{H}(K_\SSS)$
and $\Delta = H \cap \Gamma$. Let $\pi_H: H \rightarrow H/\Delta$ be
the natural projection. Remark that $\mathbf{H}$ is a reductive
group \cite[13.17, Corollary 2]{Borel}. Choose a maximal $K$-split
torus $\widetilde{\se}$ of $\bf H$. Then $\widetilde{\se} \supset
\se$ and there exists $q \in \G(K)$ such that
\begin{equation}
\label{6.1} \widetilde{\se} = q^{-1} \mathbf{D}q.
\end{equation}
Denote $\widetilde{S}_v = \widetilde{\se}(K_v)$, $v \in \SSS$, and
$\widetilde{S} = \widetilde{\se}(K_\SSS)$. There exists $h =
(h_v)_{v \in \SSS} \in H$ such that $h^{-1}_v \widetilde{S}_v h_v
\subseteq g^{-1}_v T_v g_v$ for every $v \in \SSS$. Denote
$\widetilde{T}_v =h_v g^{-1}_v T_v g_v h_v^{-1}$ and $\widetilde{T}
= \prod_{v \in \SSS}\widetilde{T}_v$. Then $\widetilde{S} \subset
\widetilde{T} \subset H$ and $\widetilde{T} \pi_H (h)$ is a closed
locally divergent orbit. Suppose for a moment that the theorem is
valid for $\mathbf{H}$. Then the conditions (i) and (ii) in the
formulation of the theorem are automatically fulfilled because
$\rank_K \G = \rank_K \mathbf{H}$ and $\rank_{K_v} \G = \rank_{K_v}
\mathbf{H}$, $v \in \SSS$. Since $h = z d$, where $z \in
\mathcal{N}_H(\widetilde{S})$ and $d \in \mathbf{H}(K)$, using
(\ref{6.1}), we obtain
\begin{align*}
D = gh^{-1} \widetilde{S}hg^{-1}  = & gd\widetilde{S}d^{-1}g^{-1} = \\
 = gdq^{-1}Dqd^{-1}&g^{-1}.
\end{align*}
Therefore, $g \in \mathcal{N}_G(D)\G(K)$, which proves (iii). The
above discussion reduces the proof to the case when $\mathbf{S}$ is
a central $K$-split torus in $\G$. In this case $\G$ is an almost
direct product over $K$ of $\mathbf{S}$ and a reductive $K$-group.
Factorizing by $\se$, we can further reduce the proof to the case
when $\mathbf{S}$ is trivial.

So, in order to complete the proof of the theorem, it is enough to
consider the case when $T_\RR \pi(g)$ is a divergent orbit. The rest
of the proof breaks in two cases according to whether or not the
assumptions in the formulation of Proposition \ref{split torus} are
satisfied.

Assume that for every $K$-subalgebra $\bb$ of $\LL$ containing
$\Lie(\mathbf{D})$ the intersection
$\mathrm{pr}_\mathcal{R}(\mathfrak{g}_x) \cap
\mathfrak{b}_{\mathcal{R}}$, where $x = \pi(g)$, contains a
horospherical subset. Then (iii) follows from Proposition \ref{split
torus}(b), and (ii) from Proposition \ref{split torus}(c) and
Theorem \ref{volumes}(b). The condition (i) follows easily from
(ii), (iii) and Proposition \ref{K-torus}.

Now assume the contrary, that is, that there exists a minimal
parabolic $K$-subalgebra $\bb$ of $\LL$ containing
$\Lie(\mathbf{D})$ and such that
$\mathrm{pr}_\mathcal{R}(\mathfrak{g}_x) \cap
\mathfrak{b}_{\mathcal{R}}$ does not contain a horospherical subset.
We will prove that this assumption leads to contradiction.  (As in
\cite{Tomanov-Weiss}, our argument is inspired by Margulis' one,
cf.\cite[Appendix]{Tomanov-Weiss}.) Let
{\mbox{\boldmath$\mathfrak{u}^-$}} be the unipotent radical of the
minimal parabolic $K$-subalgebra opposite to $\bb$. For every
positive integer $n$ we let ${B}_{n}$ be a ball of radius $n$ in
$\mathfrak{g}$. Since $\mathfrak{g}_x$ is discrete in
$\mathfrak{g}$, the family of the horospherical subsets in
$\mathrm{pr}_\mathcal{R}(\mathfrak{g}_x) \cap
\mathfrak{b}_{\mathcal{R}}$ is finite. In view of this and the
assumption that $\mathrm{pr}_\mathcal{R}(\mathfrak{g}_x) \cap
\mathfrak{b}_{\mathcal{R}}$ does not contain horospherical subsets,
for every $n$ there exists an element $s_n \in D_\RR$ such that
$\mathrm{Ad}(s_n)$ acts as an expansion on $\mathfrak{u}_{\RR}^-$
and
\begin{equation}
\label{0} \mathrm{Ad}(s_n)\mathcal{M} \nsubseteq {B}_{n}
\end{equation}
for every horospherical subset $\mathcal{M} \subset \mathfrak{g}_x
\cap {B}_{n}$.

Using Proposition \ref{intersection1}(a), we fix a compact
neighborhood $W_0$ of 0 in $\mathfrak{g}$ such that $W_0 \subset
{B}_{n}$ and for every $x \in G/\Gamma$ the subalgebra of
$\mathfrak{g}$ spanned by $\mathfrak{g}_x \cap W_0$ is unipotent.

Proposition \ref{pushing out} and the choice of $W_0$ imply that
there exist a constant $\tau
> 1$ and a finite set $t_1, \ldots, t_l$ in $D_\RR$ such that for
every $y \in G/\Gamma$ there exists $t \in \{t_1, \ldots, t_l\}$
satisfying
\begin{equation}
\label{1} \|\mathrm{Ad}(t)a\|_\RR \geq \tau\|a\|_\RR, \forall a \in
\mathfrak{g}_y \cap W_0.
\end{equation}

We put
\begin{displaymath}
W = W_0 \bigcap \big( \bigcap_{i = 1}^l \mathrm{Ad}(t_i) W_0 \big).
\end{displaymath}

Given a positive $n \in \N$, we define inductively a {\it finite}
sequence $p_0, p_1, \ldots, p_{r_n}$ as follows. We put $p_0 = s_n$.
Assume that $p_0, p_1, \ldots, p_{i}$ are already defined. If
$\mathrm{Ad}(p_{i} \ldots p_0)(\mathfrak{g}_x) \cap W$ does not
contain a horospherical subset then $p_0, p_1, \ldots, p_{i}$ is the
required sequence. If not, we put $p_{i+1} = t$, where $t$ satisfies
(\ref{1}) with $y = p_{i} \ldots p_0 x$. With the same $y$ and $p_{i
+ 1}$, remark that if $b \in \mathfrak{g}_y$ and $b \notin W_0$ then
$\mathrm{Ad}(p_{i+1})b \notin W$. This and (\ref{1}) imply the
following

$\bf{Claim}$: If $p_0, p_1, \ldots, p_{r}$ are already defined, $0
\leq i < r$, $y = p_{i} \ldots p_0 x$, $b \in \mathfrak{g}_y$ and $b
\notin W_0$ then $\mathrm{Ad}(p_j \ldots p_{i+1})b \notin W$ for
every $j$ such that $i \leq j \leq r$.

The claim implies that the cardinality of $\mathrm{Ad}(p_{i} \ldots
p_0)(\mathfrak{g}_x) \cap W$ does not increase with $i$ and,
moreover, the sequence $\{p_i\}$ is finite. Put $g_n = p_{r_n}
\ldots p_1 p_0$. It follows from Proposition \ref{intersection1}(b)
that the sequence $\{ g_n x \}$ is bounded in $G/\Gamma$. Since the
orbit $T_\RR x$ is divergent, the sequence $\{g_n\}$ is bounded in
$T_\RR$. Also note that, given the above definition of $s_n$, the
sequence $\{s_n\}$ is unbounded. Again by Proposition
\ref{intersection1}(b), passing to a subsequence,  we may assume
that $r_n > 0$ for all $n$.

Let $h_n = p_{r_n}^{-1}g_n$ and $\mathcal{M}_n$ be a horospherical
subset of $\mathrm{Ad}(h_n)(\mathfrak{g}_x) \cap W$. Assume that
$\mathrm{Ad}(h_n^{-1})(\mathcal{M}_n) \subset {B}_{n}$. Then it
follows from (\ref{0}) that $\mathrm{Ad}(p_0h_n^{-1})(\mathcal{M}_n)
\nsubseteq {B}_{n}$. The Claim implies that $\mathcal{M}_n
\nsubseteq W$, which contradicts the choice of $\mathcal{M}_n$.
Therefore,
$$
\mathrm{Ad}(h_n^{-1})(\mathcal{M}_n) \nsubseteq {B}_{n}.
$$
Since $\mathcal{M}_n \subset W$ and $W$ is compact, the sequence
$\{h_n^{-1}\}$  is not bounded. Therefore, $\{g_n\}$ is not either.
Contradiction. \qed

\subsection{Remarks} (a) It follows from the proof of Theorem \ref{thm2} that
if $\# \SSS > 1$ and the orbit $T x$ $\mathrm{(}$where $T =
T_\SSS$$\mathrm{)}$ is closed and locally divergent then the Zariski
closure of $g^{-1}T g \cap \Gamma$ in $\G$ contains a maximal
$K$-split torus.

(b) Since $\mathcal{N}_\G(\de)(K)$ meets every coset of the quotient
$\mathcal{N}_\G(\de)/\mathcal{Z}_\G(\de)$, %\cite{Borel}
we have that $\mathcal{Z}_{G}(D_v)\G(K) = \mathcal{N}_{G}(D_v)\G(K)$
for every $v$. On the other hand, it is easy to see that
$\mathcal{Z}_{G}(D_\RR)\G(K) \subsetneqq
\mathcal{N}_{G}(D_\RR)\G(K)$ whenever $\# \RR > 1$ and $\G$ is a
semisimple $K$-isotropic group.

\subsection{Proof of Theorem \ref{thm1}} Let us first prove (b).
Since the divergent orbits are locally divergent and closed we can
apply Theorem \ref{thm2}. If $\RR$ is not a singleton it follows
from Theorem \ref{thm2} (i) that $\RR = \SSS$. Also it follows from
Theorem \ref{thm2} (ii) an (iii) that $T_\SSS$ is a compact
extension of $D_\SSS$. So, $D_\SSS \pi(g)$ diverges. This
contradicts Proposition \ref{K-torus}. Therefore $\RR = \{v\}$.
Again by Theorem \ref{thm2}, $\rank_{K_v} \G = \rank_{K} \G$ and $g
\in \mathcal{N}_G(D_v)\G(K_v)$. Now in order to complete the proof
of (b) it remains to apply the remark 6.2 (b).

Let us prove (a). The implication $\Leftarrow$ follows trivially
from Propositions \ref{closed orbits1} and \ref{K-torus1}. Suppose
that $T_\RR \pi(g)$ is closed. If $T_\RR \pi(g)$ is divergent it
follows from (b) that $\RR$ is a singleton. Let $\RR = \{v\}$. Since
$\rank _{K_v} \G = \rank_K \G$, $T_v$ is a compact extension of
$D_v$. But $g = zq$, where $z \in \mathcal{Z}_G(D_v)$ and $q \in
\G(K)$. Therefore $g^{-1}T_v g$ is a compact extension of
$\mathbf{L}(K_v)$, where $\mathbf{L} = q^{-1}\de q$, which proves
(a) when $T_\RR \pi(g)$ is divergent. Let $T_\RR \pi(g)$ be not
divergent. Then $g^{-1}Tg \cap \Gamma$ is not finite, in particular,
$\RR = \SSS$. Let $\bf L$ be the connected component of the Zariski
closure of $g^{-1}Tg \cap
\Gamma$ in $\G$. Set $\mathbf{H} = \mathcal{Z}_\G(\mathbf{L})$.%, $H
%= \mathbf{H}(K_\SSS)$ and $\Delta = \mathbf{H}(\OO)$.
\ Since $\mathbf{H}$ is an almost direct product over $K$ of
$\mathbf{L}$ and of a reductive $K$-group, factorizing by
$\mathbf{L}$, we can reduced the proof to the case when $\mathbf{L}$
is trivial. In the latter case either $T_\SSS$ is compact and there
is nothing to prove or $T_\SSS \pi(g)$ is divergent. This complets
the proof of (a). \qed

\subsection{Proof of Corollaries \ref{S-div}, \ref{thm3} and \ref{thm3'}}
Corollary \ref{S-div} follows from Theorem \ref{thm1} (a) and Remark
6.2 (a), and Corollary \ref{thm3} follows from Theorem \ref{thm2}
and Remark 6.2 (b).

Let us prove Corollary \ref{thm3'}. The part (a) is immediate from
Theorem \ref{thm2}. In order to prove (b), remark that
$\big(\mathcal{N}_\G(\mathbf{D}) \times
\mathcal{N}_\G(\mathbf{D})\big) \mathrm{diag}(\G) \varsubsetneq \G
\times \G$, where $\mathrm{diag}(\G)$ is the diagonal imbedding of
$\G$ into $\G \times \G$. Therefore, there exists $(g_1, g_2) \in
(\G \times \G)(K)$ such that $(g_1, g_2) \notin
\big(\mathcal{N}_\G(\mathbf{D}) \times
\mathcal{N}_\G(\mathbf{D})\big) \mathrm{diag}(\G)$. Let $v_1$ and
$v_2$ be two different valuations in $\SSS$ and let $g = (g_v)_{v
\in \SSS} \in G$ be such that $g_{v_1} = g_1$, $g_{v_2} = g_2$ and
$g_v = 1$ for all $v \in \SSS \setminus \{v_1, v_2\}$. It follows
from Theorem \ref{thm2} (iii) and  Proposition \ref{K-torus} that
the orbit $T_\RR \pi(g)$ is locally divergent but not closed. \qed

\subsection{Remark} In connection with Corollary \ref{thm3'} (a),
note that if $G$ is a real $\Q$-algebraic group and $D_\infty$ is an
$\R$-split algebraic torus of  $G$ with $\dim D_\infty
> \rank_\Q G$, it was proved by B.Weiss \cite{Weiss1} that there are
no divergent orbits for the action of $D_\infty$ on $G/\Gamma$. The
following generalization of this result is proved \cite{Toma2}: Let
$G$ and $\Gamma$ be as in the formulation of Theorem \ref{thm1}, $v
\in \SSS$ and $\de_v$ be a $K_v$-split torus of $\G$. Assume that
$\dim \de_v > \rank_K \G$. Then $G/\Gamma$ does not admit divergent
orbits for the action of $D_v = \de_v(K_v)$.

\section{Number theoretical application}

\medskip
Let $K_\SSS[\ \vec{x}\ ]$ be the ring of polynomials in $n$
variables $\vec{x} =$ $(x_1,\ldots,x_n)$ with coefficients from the
topological ring $K_\SSS$. Let $f(\vec{x}) = l_{1}(\vec{x})\ldots
l_{m}(\vec{x})$  $\in K_\SSS[\ \vec{x}\ ]$, where
$l_{1}(\vec{x}),\ldots,$ $ l_{m}(\vec{x})$ are linearly independent
over $K_\SSS$ linear forms.

The following is a reformulation of Theorem \ref{S-inverse} from the
Introduction:

\medskip

\begin{thm}
\label{application1} With the above notation and assumptions, %Let
suppose that $f(\OO^n)$ is a discrete subset of $K_\SSS$. Then
$f(\vec{x}) = \alpha g(\vec{x})$ for some $\alpha \in K_\SSS^*$ and
some $g(\vec{x}) \in \OO[\ \vec{x}\ ]$ .
\end{thm}

\medskip

The following examples show that the hypotheses in the formulations
of Theorem \ref{application1} are essential and can not be omitted.

\medskip

$\mathbf{Examples}.$ Let $\alpha \in \R$ be a badly approximable
number, i.e. there exists a $c = c(\alpha) > 0$ such that
$$
\left|\alpha - \frac{p}{q}\right| \geq \frac{c}{q^2}
$$
for all $p/q \in \Q$. (Recall that the quadratic irrationals, such
as $\sqrt{2}$, and the golden ratio $(\sqrt{5} + 1)/2$ are badly
approximable.) Consider the form $f(x, y) = x^2(\alpha x - y)$. Then
the set of values of $f$ at the integer points is discrete {\it but}
$f$ is not a multiple of a form with rational coefficients. The
reason is that $f$ is a product of linearly {\it dependent} linear
forms.

The hypothesis that $f$ is decomposable is also essential. In order
to see this it is enough to consider a form $f(x, y) = x^2 + \beta
y^2$ where $\beta$ is a positive irrational real number. It is
obvious that $f(\Z^2)$ is discrete in $\R$.

\medskip

We put $\G = \mathbf{SL}_n$. So, $G = \mathbf{SL}_n(K_\SSS)$ and
$\Gamma = \mathbf{SL}_n(\OO)$.) The group $G$ is acting on $K_\SSS[\
\vec{x}\ ]$ according to the law $(\sigma f)(\vec{x}) =
f(\sigma^{-1}\vec{x})$, where $\sigma \in G$ and $f \in K_\SSS[\
\vec{x}\ ]$. We denote $f_0(\vec{x}) = x_1x_2...x_m$. It is clear
that if $f \in K_\SSS[\ \vec{x}\ ]$ is as in the formulation of
Theorem \ref{application1} then $f(\vec{x}) = \alpha(\sigma
f_0)(\vec{x})$ for some $\sigma \in G$ and $\alpha \in K_\SSS^*$. We
will denote by $H_f$ the stabilizer of $f$ in $G$.

\medskip

We precede the proof of Theorem \ref{application1} by the following
general proposition.

\begin{prop}
\label{integer points2} Let $f(\vec{x}) =(\sigma f_0)(\vec{x})$ for
some $\sigma \in G$ . Assume that $f(\OO^n)$ is a discrete subset of
$K_\SSS$. Then $H_f\pi(e)$ is closed in $G/\Gamma$.
\end{prop}

\medskip

{\bf Proof.} Let $\pi(a)$,  $a \in G$, belong to the closure of
$H_f\pi(e)$. Fix a sequence $h_i \in H_f$ such that $\lim_{i
\rightarrow \infty} h_i\pi(e)= \pi(a)$. There exist $\gamma_i \in
\Gamma$ and $b_i \in G$ such that $\lim_{i \rightarrow \infty} b_i
=e$ and $h_i\gamma_i = b_i a$. Since $ f(\OO^n)$ is  discrete, for
every $\vec{z} \in \OO^n$ there exists a real number $c(\vec{z}) >
0$ such that
\begin{equation}
\label{3.1}
 f(\gamma_i\vec{z}) = f(h_i\gamma_i\vec{z}) =
f(b_i a \vec{z}) = f(a\vec{z}) \in f(a \OO^n) \cap f(\OO^n)
\end{equation}
for all $i > c(\vec{z})$.

Let $\chi_1, \chi_2,...,\chi_l \in K[\ \vec{x}\ ]$ be the set of all
monomials of degree $m$. We consider $\chi_1,\chi_2,...,\chi_l$ as
homomorphisms of multiplicative groups ${K^*}^n \rightarrow K^*$.
Since $\chi_1,\chi_2,...,\chi_l$ are linearly independent over $K$,
i.e. whenever we have a relation
$$
\alpha_1\chi_1 + \alpha_2\chi_2 + \ldots + \alpha_l\chi_l = 0,
$$
with $\alpha_i \in K$ then all $\alpha_i = 0$, there exist
$\vec{z}_1, \vec{z}_2,..., \vec{z}_l \in \OO^n$ such that
$\det(\chi_k(\vec{z}_s)) \neq 0$. In view of (\ref{3.1}), there
exists $c > 0$ such that
\begin{equation}
\label{3.5} f(b_ia\vec{z}_s) = f(a\vec{z}_s)
\end{equation}
for all $s$ and $i > c$.

The form $f$ can be regarded as a collection of forms $f_v \in K_v[\
\vec{x}\ ], v \in \SSS$. Since $\det(\chi_k(\vec{z}_s)) \neq 0$,
using (\ref{3.5}), we get that
$$
f_v(b_{iv} a_v \vec{x}) = f_v(a_v\vec{x})
$$
for all $v \in \SSS$ and $i > c$, where $b_{iv}$ is the
$v$-component of $b_i$ and $a_v$ is the $v$-component of $a$. Hence
$b_i \in H_f$ for all $i > c$. So, we obtain that
$$
\pi(a) = b_i^{-1}h_i \pi(e) \in H_f \pi(e),
$$
which proves that $H_f\pi(e)$ is closed. \qed

\medskip

Given a subgroup $L$ of  $G$, we will write $L_u$ for the subgroup
generated by the Zariski closed in $G$ unipotent subgroups of $L$.

The following is a particular case of Theorem 3 from \cite{Toma}.

\begin{prop}
\label{tomanov} Let $L$ be a closed (for the Euclidean topology)
subgroup of $G$. Assume that $L\pi(e)$ is closed and $L_u\pi(e)$ is
dense in $L\pi(e)$. Let $\pe$ be the connected component of the
Zariski closure of $L \cap \Gamma$ in $\G$ and let $P =
\pe(K_\SSS)$. Then
\begin{enumerate}
\item[(i)]
$P \supset L_u$ and there exists a subgroup of finite index $P'$ in
$P$ such that $L\pi(e) = P'\pi(e)$;
\item[(ii)] If $\mathbf{Q}$ is a proper normal $K$-subgroup of $\pe$, there
exists $v \in \SSS$ such that $(\pe/\mathbf{Q})(K_v)$ contains a
unipotent element different from the identity.
\end{enumerate}
\end{prop}

\medskip

{\bf Proof of Theorem \ref{application1}. } Let $H_0$ be the Zariski
connected component of $H_{f_0}$. It is easy to see that

\begin{equation}
\label{3.2} H_{0} = \left\{\left(
\begin{array}{cc}
d&a\\
0&s\\
\end{array}
\right)|\  d \in D_m,\  a \in \mathrm{M}_{m\times (n-m)}(K_\SSS)
\text{ and } s \in \mathrm{SL}_{n-m}(K_\SSS)   \right\},
\end{equation}
where $D_m$ is the group of all diagonal matrices in
$\mathrm{SL}_m(K_\SSS)$. Since $f = \sigma f_{0}$, we have that
$$ H = \sigma H_{0}\sigma^{-1}$$
is the Zariski connected component of
$H_f$.

Let $\mathcal{F}_m$ be the $K_\SSS$-module of all homogeneous
polynomials of degree $m$ in $K_\SSS[\vec{x}]$. A simple calculation
shows that $K_\SSS f_0$ is the submodule of all $H_0$-invariant
elements in $\mathcal{F}_m$. Therefore,
\begin{equation}
\label{3.3} K_\SSS f = \{h \in \mathcal{F}_m | \sigma h = h, \forall
\sigma \in H\}.
\end{equation}

It follows from \cite[Theorem 2]{Ratner} that there exists a closed
subgroup $L$ of $G$ such that $L\pi(e) = \overline{H_u \pi(e)}$. Let
$\pe$ be the connected component of the Zariski closure of $L \cap
\Gamma$ in $\G$ and let $P = \pe(K_\SSS)$. By Proposition
\ref{tomanov}, $L\pi(e) = P'\pi(e)$ where $P'$ is a subgroup of
finite index in $P$. On the other hand, since $H_f\pi(e)$ is closed
(Proposition \ref{integer points2}) and $H$ has finite index in
$H_f$, $H \pi(e)$ is also closed. Therefore, $P' \subset H$. Since
$H_u \subset P'$, it follows from Proposition \ref{tomanov} (ii) and
from the description (\ref{3.2}) of $H_0$ that $H_u = P$ and
$L\pi(e) = P\pi(e)$.

Let $\mathbf{Q}$ be the commutator subgroup of
$\mathcal{N}_\G(\pe)$. It follows from (\ref{3.2}) that $\mathbf{Q}$
is a semidirect product over $K$ of $\pe$ and of an algebraic group
$\mathbf{R}$ defined over $K$ which is isomorphic over $K_v$ to
$\mathbf{SL}_m$ for all $v \in \SSS$. (Note that $\mathbf{R}$ is
isomorphic to $\mathbf{SL}_m$ over a finite extension of $K$ but, in
general, $\mathbf{R}$ is not isomorphic to $\mathbf{SL}_m$ over $K$
itself.) Let $R = \prod_{v \in \SSS}\mathbf{R}_v(K_v)$ and $T = R
\cap H$. Then $T = \prod_{v \in \SSS}\mathbf{T}_v(K_v)$, where
$\mathbf{T}_v$ is a maximal $K_v$-split torus in $\mathbf{R}$, and
$H = T P$. Since the projection of $H$ into $Q/(Q \cap \Gamma)$,
where $Q = \mathbf{Q}(K_\SSS)$, is closed, the projection of $T$
into $R/(R \cap \Gamma)$ is closed too. Applying Theorem \ref{thm1},
we get a torus $\mathbf{T}$ in $\mathbf{R}$ defined over $K$ such
that $T = \mathbf{T}(K_\SSS)$. Therefore, $H = \mathbf{H}(K_\SSS)$,
where $\mathbf{H} = \mathbf{T}\pe$ is an algebraic group defined
over $K$.

It follows from the above that $\mathbf{H}(K)$ is Zariski dense in
$H$. Note that given $\sigma \in \mathbf{H}(K)$ the coefficients of
all $h \in \mathcal{F}_m$ such that
$$\sigma h = h$$ can be regarded as the
space of solutions of a system of linear equations with coefficients
from $K$. Therefore, in view of (\ref{3.3}), there exist $g(\vec{x})
\in \OO[\vec{x}]$ and $\alpha \in K_\SSS^*$ such that $f(\vec{x}) =
\alpha g(\vec{x})$. \qed

\end{document}